\documentclass[11pt,margin=0.8pt]{article}
\usepackage{pkg}
\usepackage{wrapfig}
\usepackage[toc,page]{appendix}
\PassOptionsToPackage{hyphens}{url}
\usepackage[hypertexnames=false]{hyperref}
\usepackage{color}
\usepackage[dvipsnames]{xcolor}
\usepackage{latexsym}
\usepackage{dsfont}
\usepackage{appendix}
\usepackage{amssymb}
\usepackage{subcaption}
\usepackage[margin=1in]{geometry}
\usepackage{float}
\usepackage[sort&compress,numbers]{natbib}
\usepackage{amsmath, amsfonts,amssymb,euscript,array,amsfonts,mathrsfs}
\usepackage{mathabx}
\usepackage{enumerate}
\usepackage[dvipsnames]{xcolor}
\usepackage{amsthm}
\usepackage[ruled,vlined]{algorithm2e}

\usepackage[T1]{fontenc}
\usepackage[utf8]{inputenc}
\usepackage{hhline}
\usepackage{graphicx}
\usepackage{verbatim}
\usepackage{eurosym}
\usepackage{dsfont}
\usepackage{mathtools}
\usepackage{soul}

\allowdisplaybreaks

\DeclareMathOperator{\E}{\mathbb{E}}

\DeclareMathOperator{\R}{\mathbb{R}}

\DeclareMathOperator{\dd}{{\rm d}}

\makeatletter
\@addtoreset{equation}{section}

\newcommand*{\rom}[1]{\expandafter\@slowromancap\romannumeral #1@}
\makeatother
\newcommand{\vertiii}[1]{{\left\vert\kern-0.25ex\left\vert\kern-0.25ex\left\vert #1 
		\right\vert\kern-0.25ex\right\vert\kern-0.25ex\right\vert}}

\newcommand{\ba}{\begin{array}}
	\newcommand{\ea}{\end{array}}
\newcommand{\be}{\begin{equation}}
	\newcommand{\ee}{\end{equation}}
\newcommand{\eeqn}{\end{equation}}
\newcommand{\bea}{\begin{eqnarray}}
\newcommand{\eea}{\end{eqnarray}}
\newcommand{\beaa}{\begin{eqnarray*}}
\newcommand{\eeaa}{\end{eqnarray*}}

\def\a{\alpha}

\def\d{\delta}

\def\m{\mu}
\def\n{\nu}
\def\si{\sigma}

\def\f{\varphi}

\def\o{\omega}
\def\h{\widehat}

\def\Th{\Theta}

\def\O{\Omega}

\def\cA{{\cal A}}

\def\cC{{\cal C}}
\def\cD{{\cal D}}

\def\cF{{\cal F}}
\def\cG{{\cal G}}

\def\cJ{{\cal J}}

\def\cL{{\cal L}}

\def\cN{{\cal N}}

\def\cQ{{\cal Q}}

\def\cV{{\cal V}}

\def\hC{\mathbb{C}}

\def\hE{\mathbb{E}}
\def\hF{\mathbb{F}}

\def\hL{\mathbb{L}}

\def\hN{\mathbb{N}}
\def\hO{\mathbb{O}}
\def\hP{\mathbb{P}}
\def\hQ{\mathbb{Q}}
\def\hR{\mathbb{R}}

\def\hX{\mathbb{X}}

\def\sA{\mathscr{A}}
\def\sB{\mathscr{B}}

\def\sE{\mathscr{E}}

\def\sM{\mathscr{M}}

\def\sP{\mathscr{P}}
\def\sQ{\mathscr{Q}}

\def\d{\mathrm{d}}

\def\no{\noindent}

\def\ms{\medskip}
\def\bs{\bigskip}
\def\q{\quad}
\def\qq{\qquad}

\def\pa{\partial}
\def\cd{\cdot}

\def\liminf{\mathop{\underline{\rm lim}}}

\def\pa{\partial}
\def\h{\widehat}
\def\wt{\widetilde}
 
\def\cd{\cdot}

\def\1{{\bf 1}}

\def\lf[{{[\![}}
\def\Blf[{{\Big[\!\!\Big[}}
\def\ri]{{]\!]}}
\def\Bri]{{\Big]\!\!\Big]}}
\def\reff{\eqref}
\def\dis{\displaystyle}

\def\wt{\widetilde}

\newtheorem{Theorem}{Theorem}[section]
\newtheorem{Definition}[Theorem]{Definition}
\newtheorem{Proposition}[Theorem]{Proposition}
\newtheorem{Assumption}[Theorem]{Assumption}
\newtheorem{Lemma}[Theorem]{Lemma}

\newtheorem{Remark}[Theorem]{Remark}
\newtheorem{Example}[Theorem]{Example}

\newtheorem{Problem}{Problem}

\def\mi{{\mu_{\rm ini}}}
\def\mt{{\mu_{\rm tar}}}

\begin{document}
\title{Soft-constrained Schr\"odinger bridge for generative AI: numerical algorithms and approximation guarantees}

\author{Jin Ma\thanks{Department of Mathematics, University of Southern California, Los Angeles, CA, 90089. Email: jinma@usc.edu.
		This author is supported in part by US NSF grants DMS\#2205972 and \#2510403.}, ~Ying Tan\thanks{Department of Decision of Analytics and Operations,
		City University of Hong Kong, Hong Kong. Email: yingtan@cityu.edu.hk.}, Renyuan Xu\thanks{Department of Management Science \& Engineering, Stanford University, Stanford, CA 94305. Email:
		renyuanxu@stanford.edu. This author is supported in part by the NSF CAREER Award DMS\#2524465 and NSF grant DMS\#2602037.}, ~and Qinxin Yan\thanks{Department of Industrial Engineering \& Operations Research, UC Berkeley, Berkeley, CA 94709. Email: qinxinyan@berkeley.edu.}}
\date{\today}
\maketitle

\begin{abstract}
	We study a soft-constrained formulation of the Schr\"odinger bridge problem, in which the terminal distribution constraint is relaxed into a penalty on the discrepancy between the terminal law of the controlled process and a target distribution. This leads to an unconstrained McKean--Vlasov stochastic control problem, where the interaction enters only through the terminal cost. We characterize the optimizer via a McKean--Vlasov forward-backward stochastic differential equation and establish well-posedness under very mild assumptions. We also analyze two approximation regimes relevant to data-driven computation: approximation of the target law by empirical samples, and finite-particle approximation of the McKean--Vlasov control problem. For both, we derive convergence results and error estimates. Finally, building on the well-posedness and approximation results, we investigate two neural numerical methods: one based on the forward-backward system and the other on direct control optimization. Both methods are tested on state-of-the-art generative tasks. Our results provide both an algorithmic foundation and its theoretical guarentee for soft-constrained Schr\"odinger bridge problems in sample-based generative modeling.
\end{abstract}

\section{Introduction}

The Schr\"odinger bridge problem provides a principled way to construct the most likely stochastic interpolation between two probability distributions relative to a reference diffusion \citep{schrodinger1931uber,schrodinger1932sur}. Recently, it has become increasingly influential in generative AI, where one seeks to transport a prior distribution to a target data distribution through a controlled stochastic dynamics \citep{de2021diffusion,vargas2021solving,wang2021deep}. It links diffusion-based and flow-based generative modeling with stochastic control and optimal transport. At the same time, despite its conceptual elegance, the classical Schr\"odinger bridge problem remains difficult to solve numerically in large-scale and high-dimensional settings, mainly because of the hard terminal distribution constraint and the iterative procedures typically used to solve it.

In its classical stochastic control formulation, the Schr\"odinger bridge seeks an adapted drift modification that minimizes kinetic energy while steering the state process from a prescribed initial law to a prescribed terminal law. This produces a path measure that is closest to the reference diffusion in relative entropy among all path measures satisfying the marginal constraints. Enforcing both initial and terminal marginals can be computationally demanding.
For example, iterative proportional fitting schemes alternate between
updates that enforce the initial and terminal marginals, with each
update approximated numerically in practice%
~\citep{de2021diffusion,vargas2021solving,shi2024diffusion}. Moreover, the convergence rate and computational complexity of such procedures are  not well understood in high dimension. These challenges become especially important in generative AI applications, where the target distribution is accessible only through samples and the dimension of the latent space is typically large.

Motivated by these challenges, a \emph{soft-constrained} Schr\"odinger bridge formulation has been recently investigated ~\cite{garg2024soft}, ~\cite{hernandez2023propagation}  and~\cite{ma2025schr}, in which the hard terminal marginal constraint is replaced by a penalty on the discrepancy between the terminal law of the controlled process and the target distribution. More precisely, instead of forcing the terminal law to equal the target law exactly, we add to the cost function a terminal penalty functional $k\,G(\mu_T)$, where $G$ vanishes only at the target measure and $k>0$ controls the strength of the constraint. It is shown in~\cite{hernandez2023propagation} and~\cite{ma2025schr} and when $k\to\infty$, the soft-constrained Schr\"odinger problem converges back to the original Schr\"odinger problem. This relaxation converts the original problem into an unconstrained McKean--Vlasov optimal control problem, since the terminal cost now depends on the law of the terminal state. The resulting formulation is more flexible computationally, admits a broad class of penalty functionals, and is naturally suited to data-driven settings in which the target distribution is represented only through data samples. Although the relaxed formulation is promising, numerical algorithms and approximation schemes that enable large-scale implementations have not been rigorously studied in the literature. In fact, filling this gap is the primary contribution of our work.

\paragraph{Our work and contribution.}
The  soft-constrained Schr\"odinger bridge problem can also be formulated via a forward-backward system, which enjoys a special structure that could be exploited in well-posedness study as well as numerical design. By applying the stochastic maximum principle, the optimal control is represented through the adjoint process of a McKean–Vlasov forward-backward stochastic differential equation (FBSDE), where the mean-field interaction appears only in the terminal condition. Although this structure is simpler than that of a fully coupled McKean–Vlasov FBSDE, it is not covered directly by standard well-posedness results, which often rely on monotonicity assumptions or stronger regularity conditions \cite{carmona2018probabilistic}. A main theoretical objective of this paper is therefore to exploit the special structure of the soft-constrained problem to establish existence and uniqueness of the associated optimality system and thereby provide a rigorous basis for algorithm design; see Theorem~\ref{thm:global-fbsde-wellposedness}. 

A key novel ingredient of our proof is to obtain the terminal-law fixed point
from a variational problem on path space. For each candidate terminal
law, we linearize the penalty at that law and obtain an explicit
conditional Gibbs minimizer. Its terminal marginal defines a
fixed-point map. The terminal marginal of the unique minimizer of the
original path-space problem is a fixed point, which gives existence
for the FBSDE. Conversely, every FBSDE solution induces this same
minimizing path law. Uniqueness for the frozen system then gives uniqueness in law for the original FBSDE.

A second goal of the paper is to bridge the gap between theory and implementation. Even after relaxing the terminal constraint, the optimal control still depends on the terminal law through the penalty functional, and this law is typically unknown in closed form. In addition, in realistic machine-learning applications, the target distribution is observed only through finitely many samples. 
To tackle them, we study two further approximation layers. First, we analyze a data-sample approximation in which the true target law is replaced by an empirical measure built from $M$ samples, and we prove convergence of the corresponding optimal controls and values as $M\to\infty$. Second, we study a finite-particle approximation in which the McKean–Vlasov interaction is replaced by an $N$-particle system, and we obtain quantitative bounds on the induced error. These results provide theoretical justification for the particle-based numerical methods that we shall use for emerging generative AI tasks.

On the algorithmic side, we investigate two deep-learning-based methods for solving the soft-constrained problem. The first is an FBSDE-based algorithm, which parameterizes the backward process with neural networks and trains them by enforcing the terminal condition of the optimality system. The second is a direct McKean–Vlasov control algorithm, which parameterizes the feedback control itself and optimizes the  cost function through particle simulation. These two approaches correspond to two different viewpoints on the same problem: one learns the adjoint representation of the optimal feedback, while the other learns the feedback control directly from the control objective. Both methods are readily implementable in latent-space generative tasks.

\paragraph{Related literature.}
On a technical level, our formulation is closely related to \cite{hernandez2023propagation} and \cite{ma2025schr}.
The formulation of our problem is largely inspired by ~\cite{ma2025schr}, where the authors adopt a PDE and Doob's 
$h$-transform approach. Their primary goal is to establish a linear convergence rate for the soft-constrained Schr\"odinger problem toward the classical Schr\"odinger bridge problem as the penalty parameter tends to infinity. Although this linear convergence result provides useful guidance for choosing the penalty coefficient in practical applications, no specific algorithmic or numerical aspects are investigated there.
In \cite{hernandez2023propagation}, the authors use a probabilistic approach to recast a mean-field Schr{\"o}dinger bridge into a stochastic optimization problem with McKean--Vlasov dynamics, and connect the optimal control to a solution to a forward backward SDE of McKean--Vlasov type (MKV FBSDE). They consider a generalized form of the optimization problem by solving a relaxed version  and analyzing convergence. However, due to the general nature of the drift, diffusion, and running cost functions, they derived the associated MKV FBSDE without addressing the existence and uniqueness of the solution. In this work, we focus on the  soft-constrained Schr\"odinger bridge problem with specific coefficient functions motivated by applications in generative modeling, and demonstrate the solvability of the corresponding MKV FBSDE system.  Moreover, different from the techniques in \cite{hernandez2023propagation}  and~\cite{ma2025schr}, we minimize the penalized relative-entropy objective directly over path laws with a given initial marginal. The unique minimizing path law determines the
terminal-law fixed point used to establish global well-posedness of the associated McKean--Vlasov FBSDE.

From the numerical and machine learning perspective, Schr\"odinger bridge methods have attracted growing attention as a framework for data generation, distribution transport, and diffusion-based generative modeling. A major line of work studies the Schr\"odinger bridge through entropic optimal transport and iterative  Sinkhorn-type procedures \citep{nutz2022entropic,rigollet2018entropic,carlier2017convergence,jiang2026schr}, which provide a natural way to enforce the endpoint marginal constraints and have inspired several practical algorithms for simulation and generative tasks. Another line of work formulates the Schr\"odinger bridge through stochastic control, dynamic programming, or forward--backward systems \citep{dai1991stochastic,chen2016relation}, which is particularly attractive in continuous-time settings and for high-dimensional data generation. These approaches are closely connected to score-based diffusion models and controlled diffusion samplers, and provide a bridge between classical Schr\"odinger problems and modern generative modeling. For time-series generation, \citep{hamdouche2023generative} use a Schr\"odinger bridge with a path-dependent drift estimated from data.
In related work on generating stochastic paths, \citep{gao2025generativesde} use conditional diffusion models to learn the increment distributions of an unknown SDE from observed trajectories.


There are also several recent works on control-based generative modeling beyond the exact Schr\"odinger bridge setting. For instance, \cite{han2025stochastic} proposed stochastic control formulation for fine-tuning diffusion models, and prove linear convergence rate of a value-iteration algorithm based on iteratively solving the Bellman equation, which is a first result of its kind.
\cite{domingo2024adjoint} introduce adjoint matching, a memory-efficient stochastic optimal control framework for fine-tuning flow and diffusion generative models based on the maximum principle.
Reinforcement learning has also been used to fine-tune diffusion
models~\citep{Fan23,zhao2025score,han2024stochastic}. \citep{GZZ24} formulate
reward-directed training of continuous-time score-based diffusion
models as an entropy-regularized reinforcement learning problem.
More broadly, \citep{qiu2026stochastic} develops regularized mean-field control methods for global optimization over probability measures and illustrates their use in generative modeling through forward particle simulation without neural-network training.


These works share a stochastic control perspective on generative modeling. For the soft-constrained Schr\"odinger bridge, the dependence of the terminal penalty on the state distribution leads us to numerical methods for McKean--Vlasov control and the associated FBSDEs. In particular, our two numerical schemes, one based on the forward--backward representation and one based on direct optimization of the McKean--Vlasov control problem, are designed to reflect the specific mean-field terminal interaction of the soft-constrained formulation. The former numerical scheme is also analyzed in~\cite{carmona2022convergence} for more general McKean--Vlasov control and mean-field game problem. Numerical algorithms for general McKean--Vlasov control problem are also studied in \citep{soner2025learning, PW1,PW2,PW_wasserstein, DLZ}. In this work, we mainly adopt the algorithm proposed in~\cite{soner2025learning}. We complement these studies by establishing global well-posedness of the bridge FBSDE under our standing assumptions, including convexity of the terminal penalty along linear mixtures, and by deriving error bounds for approximating the target law from data samples and the controlled terminal law with finitely many particles.

\paragraph{Notations.} Throughout this paper,  we consider  a  generic Euclidean space $\hX$, and regardless of its dimension,  we denote $(\cd,\cd)$ and $|\cd|$ be its inner product and norm, respectively. We denote $\hC([0,T];\hX)$ to be the space of $\hX$-valued continuous functions defined on $[0,T]$ with the usual sup-norm, and in particular, we denote $\hC^d_T:=\hC([0,T];\hR^d)$, and  let $\sB(\hC_T^d)$ be its topological Borel field.
We shall consider the following 
{\it canonical probabilistic space}: $(\O,\cF, \hP)$, where $(\O, \cF):=(\hC^d_T, \sB(\hC^d_T))$ and $\hP\in\sP(\O)$, the space of all probability measures defined on $(\O, \cF)$. Finally, we let $\hP^0\in \sP(\O)$ be
the Wiener measure on $(\O, \cF)$;  $W_t(\o):=\o(t)$, $\o\in\O$, the canonical process; and  $\hF^0:=\{\cF^0_t\}_{t\in[0,T]}$,  where $\cF^0_t:=\sB_t(\hC_T^d):=\si\{\o(\cd\wedge t):\o\in\hC^d_T\}$, $t\in[0,T]$. We assume that $\hF^0$ has the usual $\hP^0$-augmentation so that it satisfies the {\it usual hypotheses} (cf. e.g., \cite{protter2005stochastic}), and for $p\ge 0$, and $\cG\subseteq \cF_T$, we denote $\hL^p_{\cG}(\O;\hX)$ to be all $\hX$-valued, $p$-integrable, $\cG$-measurable random variables, and $\hL^p_{\hF^0}([0,T];\hX)$ to be all $\hX$-valued, $p$-integrable, $\hF^0$-adapted processes. Finally, we denote $\sM(\hR^d)$ to be all $\si$-finite measures on $\hR^d$ and $\sP_p(\hR^d)$ to be all probability measures with finite $p$-th moment on $\hR^d$ equipped with $p$-Wasserstein metric, denoted by $W_p(\cd, \cd)$.  We shall denote $\sP(\hR^d)=\sP_1(\hR^d)$ for simplicity.

\section{Preliminary}\label{sec:preliminary}

{\color{black} We begin by recalling the classical Schr\"{o}dinger bridge problem (SBP). For any given $\mi, \mt\in\sP(\hR^d)$ and $\hP\in\sP(\O)$ an SBP is to find the solution to the following optimization problem:
	\bea
	\label{SP0}
	V(\hP):= \inf_{\hQ\in\sP(\mi, \mt)} D_{\rm KL}(\hQ\|\hP),
	\eea
	where 
	$D_{\rm KL}(\,\cd\,\|\,\cd\,)$ is the {\it Kullback-Leibler divergence} or the {\it Relative Entropy} (cf. \cite{KL}),
	defined by
	\beaa
	D_{\rm KL}(\hQ\|\hP)
	:=
	\begin{cases}
		\hE^{\hQ}\big[\log\big(\frac{\d\hQ}{\d\hP}\big)\big], & {\text{if }} \d\hQ \ll d\hP;\\
		\infty, &{\text{otherwise}},
	\end{cases}
	\eeaa
	and the ``admissible set'' $\sP(\mi, \mt)\subseteq \sP(\O)$ in (\ref{SP0}) is defined by 
	\beaa
	\sP(\mi, \mt):=\{\hQ\in\sP(\O): ~ D_{\mathrm{KL}}(\hQ\Vert \hP)<\infty, ~
	{\hQ}_{{X_0}}=\mi, ~{\hQ}_{{X_T}}  =  \mt
	\},
	\eeaa
	where $X$ is the canonical process defined on $(\hC^d_T, \sB(\hC^d_T))$. Define the relaxed set as 
	\begin{equation*}
		\sQ_{\mu_{\mathrm{ini}}} :=
		\left\{
		\hQ\in\sP(\Omega): ~ D_{\mathrm{KL}}(\hQ\Vert \hP)<\infty, ~\hQ_{X_0}=\mu_{\mathrm{ini}}\right\}\supseteq \sP(\mi, \mt).
		\label{eq:finite-entropy-path-laws}
	\end{equation*}
	
	It is well-known that  (see, e.g., \cite[\S 4.4]{chen2021stochastic},  \cite[\S 1]{chow2022dynamical},  \cite[p36]{leonard2013survey}, and \cite{ma2025schr}) that the SBP (\ref{SP0}) can be recast as the following stochastic control problem:}
\bea
\label{SB-objective}
V(\mi,\mt)=\inf_{\alpha \in \mathcal{A}}J(\alpha) = \inf_{\alpha \in \mathcal{A}}\hE^\hQ \Big[ \frac{1}{2} \int_0^T |\alpha_s|^2 \d s  \Big],
\eea
where $\hQ\in\sP(\O)$ is such that $\frac{\d\hQ}{\d\hP}=\sE(\a)$ for some $\a\in \mathcal{A}\subseteq \hL^2_{\hF^0}([0,T];\hR^d)$, and under $\hQ$ the controlled dynamics takes the form
\bea
\label{SB-dynamics}
\d X^\alpha_t = [b(t,X^\alpha
_t) +\alpha_t]\d t + \si(t)\d\wt{W}_t, \quad \hQ_{X_0^\alpha} =\mi, \q \hQ_{X_T^\alpha}= \mt,
\eea
where $\wt{W}$ is a $\hQ$-Brownian motion. {\color{black} Thus, the optimizer of \eqref{SB-objective}--\eqref{SB-dynamics} yields a diffusion interpolation between $\mi$ and $\mt$ whose path law minimizes the KL-divergence
	relative to $\hP$, the law of the reference diffusion
	\bea
	\label{P0}
	\d X_t=b(t,X_t)\,\d t+\si(t)\d B_t,\qquad \hP_{X_0}=\mi,
	\eea
	where $B$ is a $\hP$-Brownian motion.

	The main difficulty for solving such a control problem, even numerically, is to handle the ``fixed-point'' constraints at both the initial and terminal times simultaneously, albeit in the sense of distributions.
	
	Following \citep{hernandez2023propagation} and \citep{ma2025schr},
	we replace the terminal marginal constraint in~\eqref{SB-dynamics}
	by a penalty on the discrepancy between the controlled terminal law
	and $\mt$. The resulting McKean--Vlasov control problem keeps the
	prescribed initial law and has no hard terminal constraint.}

\medskip

More specifically,  given $\mt\in\sP(\hR^d)$, we first introduce the  notion of {\it penalty function}.
\begin{Definition}
	\label{penalty}
	A function $G(\cd)=G(\cdot;\mt):\sP_2(\hR^d) \mapsto [0,\infty)$ is called a (continuous) penalty function associated to $\mt\in \sP_2(\hR^d)$ if:  (1) $G(\cd)$ is continuous on $\sP_2(\hR^d)$; and (2) 
	$G(\mu;\mt) = 0$ if and only if $\mu=\mt$.
	\qed
\end{Definition}

We recall (see e.g., \cite{carmona2018probabilistic}) that, given a function $G(\cd)$, a jointly Borel version of the scalar linear functional derivative of $G$, denoted by
$ \frac{\delta G}{\delta\mu}(\m)(x)$, $(\m,x)\in\sP_2(\hR^d)\times \hR^d$, is a jointly Borel measurable function on $(\m,x)$, such that the folloing identity holds:
\begin{equation}
	\label{Deriv}
	\lim_{\varepsilon\downarrow0}  \frac{ G((1-\varepsilon)\mu+\varepsilon\nu)-G(\mu)}{\varepsilon}
	= \int_{\mathbb R^d} \frac{\delta G}{\delta\mu}(\mu)(x)(\nu-\mu)(\d x), \q \mu,\nu\in\sP_2(\mathbb R^d).
\end{equation}
Moreover, the linear {\it Bregman remainder} is defined by
\bea
\label{eq:linear-bregman-remainder}
\mathcal B_G(\nu\mid\mu) := G(\nu)-G(\mu) - \int_{\mathbb R^d}\frac{\delta G}{\delta\mu}(\mu)(x)(\nu-\mu)(\d x).
\eea
We should note that the Gateau-type scalar derivative is different from the well-known vector-valued Lions derivative $\partial_\m G(\m)$. In fact, it holds that, for each $\m\in\sP^2(\hR^d)$, 
\bea
\label{relation}
\partial_\mu G(\mu)(x)=\nabla_x\frac{\delta G}{\delta\mu}(\mu)(x), \qq x\in\hR^d.
\eea

We now consider the following family of McKean--Vlasov-type stochastic control problems.

\begin{Problem}[Soft-constrained Schr\"odinger bridge problem (SCSBP)] 
	\label{softSBP}
	For $k\in\hN$, find $\h\a^k\in\cA$ such that $V^k: =J^k(\widehat\a^k)=\inf_{\alpha \in \cA} J^k(\alpha) $, where
	\bea
	\label{RSB_objective}
	J^k(\alpha)=\hE^\hQ \left[ \frac{1}{2} \int_0^T |\alpha_s|^2 \d s  + k G({\hQ_{X^\alpha_T}})\right],
	\eea
	and   $G(\cd)=G(\cdot;\mt)$ is the given penalty  function satisfying Definition \ref{penalty}  and $\hQ\sim \a$.
	\qed
\end{Problem}

In light of \reff{SB-objective}-\reff{SB-dynamics}, we can characterize Problem~\ref{softSBP} in the following 
equivalent formulation, which will facilitate our analysis below.

For every $\hQ\in\sQ_{\mi}$, let $\alpha^{\hQ}$ denote its canonical
finite-entropy drift, namely the progressively measurable process, such that $\frac{\d\hQ}{\d\hP}\Big|_{\cF_T}=\sE(\a)=\exp\big\{\int_0^T \a^\hQ_s\d W_s-\frac12\int_0^T|\a^\hQ_s|^2\d s\big\}$, and
there exists a $\hQ$-Brownian motion $W^{\hQ}$ satisfying
\begin{equation}
	\label{eq:controlled-dynamics-follmer}
	\d X_t =\left(
	b(t,X_t)+\sigma(t)\alpha_t^{\hQ}
	\right)\d t + \sigma(t)\,\d W_t^{\hQ}, \qq\hQ_{X_0}=\mi.
\end{equation}
We shall refer to the pair $(\hQ, \a^\hQ)$ as the {\it weak admissible control}, and often write it as $\hQ\sim\alpha^{\hQ}$. Furthermore, we denote the set of weak admissible control by
\begin{equation}
	\sA_{\mathrm w} := \big\{
	(\hQ,\alpha^{\hQ}):
	\hQ\in\sQ_{\mi}
	\big\}.
	\label{eq:weak-admissible-control-class}
\end{equation}
In particular, we shall always denote the $\hP\in\sP(\O)$ be such that $\hP$ is the law of the reference diffusion (\ref{P0}) with $\hP_{X_0}=\mi$, that is $\hP\sim 0$ (or $(\hP,0)\in\sA_{\mathrm w}$). 


{\color{magenta}
}

\begin{Problem}[Path-law formulation]
	\label{Weak_SCSP}
	The weak path-law formulation of the SCSBP is
	\begin{equation}
		V^k
		:=
		\inf_{\hQ\in\sP(\O),\\\hQ_{X_0} = \mi}J^k(\hQ),
		\qquad
		J^k(\hQ)
		:=
		D_{\mathrm{KL}}(\hQ\Vert\hP)+kG(\hQ_{X_T}).
		\label{eq:path-law-objective}
	\end{equation}
\end{Problem}

\ms

We first show the equivalence between Problem~\ref{softSBP} and Problem~\ref{Weak_SCSP}.
\begin{Proposition}[Equivalence of the weak control and path-law formulations]
	\label{prop:control-path-law-equivalence}
	Let $\hP\in\sP(\O)$ be such that $\hP\sim 0$ (or $(\hP,0)\in\sA_{\mathrm w}$). 
	Then, for every $k>0$,
	\begin{equation}
		\begin{aligned}
			&\inf_{(\hQ,\alpha)\in\cA_{\mathrm w}}
			\left\{
			\frac12
			\mathbb E^{\hQ}
			\left[
			\int_0^T|\alpha_t|^2\,\d t
			\right]
			+
			kG(\hQ_{X_T})
			\right\}
			=
			\inf_{\substack{\hQ\in\sP(\O)\\ \hQ_{X_0}=\mi}}
			\left\{
			D_{\mathrm{KL}}(\hQ\Vert\hP)
			+
			kG(\hQ_{X_T})
			\right\}.
		\end{aligned}
		\label{eq:control-path-law-value-equivalence}
	\end{equation}
	In particular, the two problems have the same value, and the the optimizers correspond as follows:
	\begin{enumerate}
		\item  If $(\widehat{\hQ}^{\,k},\widehat\alpha^{\,k})$ is optimal Problem 1, then 
		$\widehat{\hQ}^{\,k}$ is optimal for Problem 2. 
		
		\item If $\widehat{\hQ}^{\,k}$ is optimal for Problem 2,  then its canonical drift $\alpha^{\widehat{\hQ}^{\,k}}$ is optimal
		Problem 1.
		
		\item  If the optimal path law for Problem 2 is unique, then the optimal canonical
		control  for Problem 1 is unique, $\d t\otimes \d\widehat{\hQ}^{\,k}$-a.e.
	\end{enumerate}
\end{Proposition}

\begin{proof}
	We first establish the finite-entropy drift representation and the corresponding entropy identity. To this end, let $\hQ\in\sQ_{\mi}$, and $\hP$ be the law of the reference diffusion (\ref{P0}). We shall denote the corresponding $\hP$-Brownian motion by $W^\hP$ for notational clarity. 
	
	Since $D_{KL}(\hQ\| \hP)<\infty$, we define $\frac{\d\hQ}{\d\hP}\Big|_{\cF_t}:=L_t$, $t\in[0,T]$, then $\{L_t\}$ is a $\hP$-martingale with $L_0=1$. Moreover, by Martingale Representation Theorem, there exists $H\in L^2_\hF([0,T])$, such that
	$ L_t=1+\int_0^t H_s\d W^{\hP}_s
	$ for every $0\le t\le T$.
	Note that, for each $t\in[0,T]$, we have $\hQ\{L_t=0\}=\hE^\hP\{L_t\1_{\{L_t=0\}}\}=0$, that is, $\hQ\{L_t>0\}=1$. Let us define $\a^{\hQ}_t=H_t/L_t\1_{\{L_t>0\}}=H_t/L_t$, 
	$t\in[0,T]$, $\hQ$-a.s. Then, under $\hP$, it holds that
	\begin{equation}
		\d L_t = L_t\alpha_t^{\hQ}\d W_t^{\hP}, \qq L_0=1.
		\label{eq:density-stochastic-logarithm}
	\end{equation}
	Thus $L$ can be written as Dol\'eans-Dade stochastic exponential $L_T=\sE(\a)=\exp\big\{\int_0^T \a^\hQ_sdW^\hP_s-\frac12\int_0^T|\a^\hQ_s|^2ds\big\}$, and consequently, 
	%
	\begin{equation}
		\mathbb E^{\hP}
		\left[L_T\log L_T\right]=\frac12 \mathbb E^{\hP}\Big[ \int_0^T L_t|\alpha_t^{\hQ}|^2\,\d t\Big]=\frac12
		\mathbb E^{\hQ}\Big[\int_0^T|\alpha_t^{\hQ}|^2\,\d t \Big]=D_{\mathrm{KL}}(\hQ\Vert\hP).
		\label{eq:entropy-energy-identity-equivalence}
	\end{equation}
	
	Furthermore, define $W_t^{\hQ}:= W_t^{\hP} - \int_0^t\alpha_s^{\hQ}\,ds$. Then, by Girsanov Theorem, the $W^\hQ$ is a $\hQ$-Brownina motion and under $\hQ$,  the reference dynamic (\ref{P0}) now reads
	\begin{equation}
		\d X_t = \left(
		b(t,X_t)+\sigma(t)\alpha_t^{\hQ}
		\right)\d t
		+
		\sigma(t)\,\d W_t^{\hQ}.
	\end{equation}
	Thus $\alpha^{\hQ}$ is an admissible weak control associated with $\hQ$. The uniqueness of $\a^\hQ$ clearly follows from that of the Radon-Nidodym derivative $L$, and Martingale Representation kernel $H$. 
	
	
	We now compare the two objective functionals. To distinguish the cost functionals in two Problems, for each $\hQ\in\sQ_{\mi}$ we denote
	$J^k(\a^\hQ)$ in \reff{RSB_objective}  by $J_{\mathrm{control}}^k (\hQ,\alpha^{\hQ})$, and that in \reff{eq:path-law-objective}  by  $J_{\mathrm{law}}^k(\hQ)$. Then it follows from \eqref{eq:entropy-energy-identity-equivalence} that 
	\begin{align}
		J_{\mathrm{control}}^k
		\left(
		\hQ,\alpha^{\hQ}
		\right)
		:=
		\frac12
		\mathbb E^{\hQ}
		\left[
		\int_0^T
		|\alpha_t^{\hQ}|^2\,\d t
		\right]
		+
		kG(\hQ_{X_T})
		=
		D_{\mathrm{KL}}(\hQ\Vert\hP)
		+
		kG(\hQ_{X_T})=
		J_{\mathrm{law}}^k(\hQ).
		\label{eq:pointwise-objective-equivalence}
	\end{align}
	
	Consequently,  the following identities yield \eqref{eq:control-path-law-value-equivalence}:
	\begin{align}
		\label{infJk}
		\inf_{\substack{\hQ\in\sP(\O)\\ \hQ_{X_0}=\mi}}
		J_{\mathrm{law}}^k(\hQ)
		&=
		\inf_{\hQ\in\cQ_{\mi}}
		J_{\mathrm{law}}^k(\hQ)
		=
		\inf_{\hQ\in\cQ_{\mi}}
		J_{\mathrm{control}}^k
		\left(
		\hQ,\alpha^{\hQ}
		\right)
		=
		\inf_{(\hQ,\alpha)\in\sA_{\mathrm w}}
		J_{\mathrm{control}}^k(\hQ,\alpha).
	\end{align}
	
	Combining \reff{eq:pointwise-objective-equivalence} and \reff{infJk} we see that the minimizers of the Problem 1 and Problem 2 must be identical,
	proving the proposition. 
\end{proof}

Note that the well-posedness of Problem \ref{softSBP} can be studied either from the stochastic control/PDE perspective or from the FBSDE perspective, which will form the basis for our numerical algorithm design. The stochastic control/PDE side is relatively well understood under various settings; see \cite{soner2025learning,ma2025schr}. However, from the FBSDE perspective, much of the existing existence-and-uniqueness theory relies on very strong model assumptions and does not fully exploit the structure of Problem \ref{softSBP}, where the mean-field interaction appears {\it only} through the terminal cost $G$.

In the next two sections, we shall establish well-posedness for the associated coupled FBSDE system, and then discuss our algorithmic approach based on both McKean--Vlasov control and FBSDE methods.

\section{Stochastic Control Problem and associated FBSDE}


Throughout this paper, we shall make the following {\it Standing Assumption} on the coefficients $b$, $\si$, as well as the penalty function $G$.

{\color{black}
	\begin{Assumption}
		\label{assump1}
		(1)  The coefficients $b:[0,T]\times\mathbb{R}^d\to\mathbb{R}^d$ and 
		$\sigma:[0,T]\to\mathbb{R}^{d\times d}$ are continuous functions, such that 
		
		\begin{itemize}
			\item[(i)] For fixed $t\in[0,T]$, the function $b\in \cC^{1,2}([0,T]\times \hR^d)$, and there exist constants $L, L' > 0$ such that, for all 
			$t\in[0,T]$, it holds that
			$$\|\nabla_x b(t,x)\|\le L, \q \|\nabla_x b(t,x)-\nabla_x b(t,y)\|\le L' |x-y|, \qq x,y\in\mathbb{R}^d.
			$$
			
			\item[(ii)] There exist $0<c<C$ such that 
			$c|x|^2\le x^T\si(t)\si^T(t)x\le C|x|^2$, $x\in\hR^d$, uniformly in $t\in [0,T]$;
		\end{itemize}
		
		\ms
		
		\no(2) 
		The penalty function $G(\cd)$ defined by Definition \ref{penalty} satisfies, for some constant $c>0$, 
		\begin{itemize}
			\item[(i)]  $G(\hP_X) \leq  c(1+\hE^\hP[|X|^2])$, $X\in \hL^2(\O;\hP)$;
			
			\item[(ii)] 	$|G(\mu)-G(\mu')|\leq c W_2(\mu,\mu')$, $\m,\m'\in \sP_2(\hR^d)$.
			
			\item[(iii)]	
			$
			G((1-\theta)\mu+\theta\nu)
			\leq
			(1-\theta)G(\mu)+\theta G(\nu),
			\qquad
			\theta\in[0,1],\quad\mu,\nu\in \P(\hR^d).
			$
			\item[(iv)] $G$ admits a scalar linear functional derivative
			$\frac{\delta G}{\delta\mu}$. After fixing an $x$-independent normalization, the
			map
			$ 
			(\mu,x)\longmapsto
			\frac{\delta G}{\delta\mu}(\mu)(x)
			$ 
			is jointly Borel measurable, $\frac{\delta G}{\delta\mu}(\mu)(\cdot)\in \mathcal C^3(\mathbb R^d)$ for every
			$\mu$, and there exists $L_G>0$ such that
			$ 
			\sum_{j=0}^3\left\|D_x^j\frac{\delta G}{\delta\mu}(\mu)\right\|_\infty
			\leq L_G
			$ for every $\mu\in\sP_2(\mathbb R^d)$.
		\end{itemize}
	\end{Assumption}

	\begin{Remark}
		\label{rem:linear-versus-lifted-convexity}
		{\rm
			Assumption~\ref{assump1}-(2)-(iii) can be easily interpreted in terms of  $\frac{\d G}{\d \m}$ (see \reff{Deriv}) by
			\begin{equation}
				G(\nu)-G(\mu)
				\geq
				\int_{\mathbb R^d} \frac{\delta G}{\delta\mu}(\mu)(x)(\nu-\mu)(\d x),\qquad \forall\,\mu,\nu\in \sP_2(\hR^d).
				\label{eq:linear-subgradient-assumption}
			\end{equation}
			Or equivalently, in terms of  Bregman's remainder (see \reff{eq:linear-bregman-remainder}):
			$ \mathcal B_G(\nu\mid\mu)\geq0$, $ \mu,\nu\in\sP_2(\mathbb R^d)$.
			\qed}
	\end{Remark}

	\paragraph{Convexity along linear mixtures and convexity of the Lions lift.}
	Assumption~\ref{assump1} requires convexity along linear mixtures of
	probability measures. In the probabilistic approach to
	McKean--Vlasov control, a common alternative is convexity of the
	Lions lift on $L^2$~\citep{carmona2018probabilistic}. Convexity along linear mixtures suffices for the convexity step in our well-posedness proof and accommodates popular penalties for generative AI tasks, such as Gaussian maximum mean discrepancy and the squared Fourier--Wasserstein metric, whose lifts need not be convex.
	
	Let us first recall that for any $\m\in\sP(\O)$, the Lions lift for function $G(\cd)$  is defined by 
	a function $\widetilde G:L^2(\widetilde{\Omega}; \hR^d)\mapsto \hR$,
	where $(\widetilde{\Omega},\widetilde{\mathcal F},\widetilde{\mathbb P})$ is some atomless probability space, such that for any $\xi\in L^2(\widetilde{\Omega};\hR^d)$ with $\m=\cL(\xi)$, it holds that $\widetilde G(\xi)= G(\mathcal L(\xi))$. 
	The convexity in the sense of \cite{carmona2018probabilistic} means the function $\widetilde G(\cd)$ is convex. When $G$ is Lions differentiable,
	the convexity of $G$ implies that 
	\bea
	\label{convexLions}
	G\bigl(\mathcal L(\xi')\bigr) - G\bigl(\mathcal L(\xi)\bigr) \geq
	\widetilde{\mathbb E}\big[
	\partial_\mu G\bigl(\mathcal L(\xi)\bigr)(\xi)
	\cdot(\xi'-\xi) \big],
	\eea
	for every pair $\xi,\xi'\in L^2(\wt{\O};\hR^d)$. We note that this inequality must
	hold for every coupling of $(\xi, \xi')$ with the marginal laws 
	$\mathcal L(\xi)$ and $\mathcal L(\xi')$, or equivalently,
	\beaa
	G\bigl(\mathcal L((1-\theta)\xi+\theta\xi')\bigr) \leq 
	(1-\theta)G\bigl(\mathcal L(\xi)\bigr)
	+\theta G\bigl(\mathcal L(\xi')\bigr),
	\qquad
	\theta\in[0,1],
	\eeaa
	and the inequality holds for all possible   joint  distribution of  $(\xi,\xi')$.
	
	By contrast, the convexity used in this paper is  along the
	linear mixtures of probability measures:
	$$
	G\bigl((1-\theta)\mu+\theta\nu\bigr)
	\leq
	(1-\theta)G(\mu)+\theta G(\nu),
	\qquad
	\mu,\nu\in\sP_2(\mathbb R^d),
	\quad
	\theta\in[0,1], 
	$$
	and the inequality depends only on the marginals $\m=\cL(\xi)$ and $\n=\cL(\xi')$. In fact, denoting 
	$
	\phi_\mu(x):=\frac{\delta G}{\delta\mu}(\mu)(x)$ for notational simplicity, we see that this leads to the subgradient inequality
	\bea
	\label{convexlinear}
	G(\nu)-G(\mu) \geq \int_{\mathbb R^d} \phi_\mu(x)(\nu-\mu)(\d x).
	\eea
	
	We remark that the two convexity conditions concern different interpolations and mutually non-inclusive. Roughly speaking, 
	the measure $(1-\theta)\mu+\theta\nu$ is simply a linear mixture of two marginals, whereas the law
	$ 
	\mathcal L((1-\theta)\xi+\theta\xi')
	$ 
	depends on the entire coupling $\mathcal L(\xi,\xi')$. They are not equal in general. 
	%
	For example, consider the variance functional for a given $\m\in\sP_2(\hR^d)$, 
	$$
	\mathcal V(\mu)
	:=
	\int_{\mathbb R^d}|x-m_\mu|^2\mu(\d x),
	\qquad
	m_\mu:=\int_{\mathbb R^d}x\mu(\d x).
	$$
	It is easy to check that  Lions lift for $\cV(\cd)$ is
	$\widetilde{\mathcal V}(\xi)=\widetilde{\mathbb E}\big[|\xi-\widetilde{\mathbb E}[\xi]|^2\big]$, where $\xi\sim \m$, and as the squared norm of a linear projection on $L^2$, it is  convex in terms of $\xi$. However, 
	\beaa
	\mathcal V\bigl((1-\theta)\mu+\theta\nu\bigr)
	=(1-\theta)\mathcal V(\mu)+\theta\mathcal V(\nu) +\theta(1-\theta)|m_\mu-m_\nu|^2\geq (1-\theta)\mathcal V(\mu)+\theta\mathcal V(\nu),
	\eeaa
	that is,  $\mathcal V(\cd)$ is generally concave, or at most linear when 
	$m_\m=m_\n$. In other words, the lifted convexity does not imply the linear convexity.

	Another important implication of the convexity is the so-called ``monotonicity conditions'' often used in the direct comparison or
	continuation arguments for McKean--Vlasov FBSDEs. 
	More precisely, lifted convexity implies $\widetilde{\mathbb E}\left[
	\left(
	\partial_\mu G\bigl(\mathcal L(\xi')\bigr)(\xi')
	-
	\partial_\mu G\bigl(\mathcal L(\xi)\bigr)(\xi)
	\right)
	\cdot(\xi'-\xi)
	\right]
	\geq 0$
	for every coupling $(\xi,\xi')$, which is often called the {\it terminal monotonicity
		condition}, but the linear convexity
	only yields the {\it scalar monotonicity condition} $\int_{\mathbb R^d}
	\left(
	\phi_\nu(x)-\phi_\mu(x)
	\right)
	(\nu-\mu)(\d x)
	\geq 0$.
	
	\medskip 
	Assumption~\ref{assump1}-(1) imposes regularity and uniform
	ellipticity conditions on the reference diffusion%
	~\citep{protter2005stochastic,zhang2017backward}.
	The next examples, inspired by~\cite[Example~2.2]{BLPR}
	and~\cite{soner2025learning}, give penalty functionals that satisfy
	Assumption~\ref{assump1}-(2) and can fail to satisfy lifted convexity.

	\begin{Example}[Gaussian max-mean discrepancy (MMD)]\label{ex:mmd-gaussian}
		Fix $\eta>0$, consider the Gaussian kernel
		$k_\eta(x,y):=\exp\big\{-\frac{|x-y|^2}{2\eta^2}\big\}$,  $ x,y\in\hR^d$.
		For $\mu\in\sP_2(\hR^d)$, define
		\begin{align*}
			G(\mu)
			:= \iint_{\hR^d\times\hR^d} k_\eta(x,x')\,\mu(\d x)\mu(\d x')
			-2\iint_{\hR^d\times\hR^d} k_\eta(x,y)\,\mu(\d x)\mt(\d y) +\iint_{\hR^d\times\hR^d} k_\eta(y,y')\,\mt(\d y)\mt(\d y').
		\end{align*}
		Since $k_\eta$ is positive definite, we have $G(\mu)\ge 0$.
		Moreover, since the Gaussian kernel is characteristic, it holds that $G(\mu)=0$ if and only if $\mu=\mt$.
		Therefore $G$ satisfies Definition~\ref{penalty}.
		
		To see that $G(\cd)$ satisfies  Assumption~\ref{assump1}-(2), we note that\\
		(i)
		$0\le G(\mu)\le 4$ for all $\mu$, hence $G(\hP_X)\le 4 \le 4(1+\hE^\hP[|X|^2])$ for all $X\in \hL^2(\O;\hP)$; \\
		(ii)  $\nabla_x k_\eta(x,y)= -\big(\frac{x-y}{\eta^2}\big)k_\eta(x,y)$, and consequently,  
		\begin{center}
			$L_\eta:=\sup_{x,y}|\nabla_x k_\eta(x,y)|=\sup_{x,y} \frac{|x-y|}{\eta^2}e^{-\frac{|x-y|^2}{2\eta^2}}
			\le \frac{1}{\eta}.$
		\end{center}
		
		Now, for $\mu,\mu'\in\sP_2(\hR^d)$, let $(X,\bar X)$ be any coupling with ${\rm Law}(X)=\mu$, ${\rm Law}(\bar X)=\mu'$,
		and let $(X',\bar X')$ be an independent copy of $(X,\bar X)$. Let $Y,Y'$ be i.i.d.\ with law $\mt$,
		independent of everything else.  Using the Lipschitz continuity of $k_\eta$ in each variable with constant $L_\eta$,
		\begin{align*}
			|G(\mu)-G(\mu')|
			&\le \E\big|k_\eta(X,X')-k_\eta(\bar X,\bar X')\big|
			+2\,\E\big|k_\eta(X,Y)-k_\eta(\bar X,Y)\big|\\
			&\le 2L_\eta\,\E|X-\bar X|+2L_\eta\,\E|X-\bar X|
			=4L_\eta\,\E|X-\bar X|.
		\end{align*}
		Taking the infimum over couplings gives $|G(\mu)-G(\mu')|\le 4L_\eta W_1(\mu,\mu')\le 4L_\eta W_2(\mu,\mu')$, proving Assumption \ref{assump1}-(2).
		
		Finally, we should note that $G$ is $L$-differentiable and its Lions derivative is given by
		\begin{center}
			$\partial_\mu G(\mu)(x)
			\,=\,2\int_{\R^d}\nabla_x k_\eta(x,x')\,\mu(\d x')
			-2\int_{\R^d}\nabla_x k_\eta(x,y)\,\mt(\d y),
			\qquad x\in\R^d.$ \hfill \qed
		\end{center}
	\end{Example}


	\begin{Example}[Fourier Wasserstein]
		\label{ex: fourier-wasserstein}
		In this example we consider the Fourier basis: 
		$e(x,\xi):=(2\pi)^{-\frac{d}{2}}e^{i\xi\cdot x}$, $x,  \xi\in \hR^d$, 
		where $i:=\sqrt{-1}$. For every $\mu\in \sP_2(\hR^d)$, denote its Fourier transform as
		\begin{equation*}
			\mathfrak{F}(\mu)(\xi):=\int_{\hR^d} e^*(x,\xi)\mu(\d x),\qquad \xi\in \hR^d,
		\end{equation*}
		where  $z^*$  denotes the complex conjugate of $z\in \hC^d$.  Clearly,  $\|\mathfrak{F}(\mu)(\cdot)\|_{\infty}\le 1$ for every $\mu\in \sP_2(\hR^d)$. Furthermore, for   $\lambda>\lceil\tfrac{d}{2}\rceil$, we define the so-called $\lambda$-th Sobolev-Wasserstein metric on $\sP_2(\hR^d)$ by
		\begin{equation}
			\rho_\lambda(\mu,\nu)^2:=\int_{\hR^d}(1+|\xi|^2)^{-\lambda}|\mathfrak{F}(\mu)(\xi)-\mathfrak{F}(\nu)(\xi)|^2 \d\xi,\qquad \mu,\nu\in \sP_2(\hR^d).
		\end{equation}
		Then, it is known that
		if $\lambda\ge \lceil\tfrac{d}{2}\rceil+1$, there exists a constant $C_d>0$, depending only on the dimension $d$, such that $\rho_\lambda(\mu,\nu)\le c_d W_1(\mu,\nu)$ (see, e.g.,  \cite[Section 2]{soner2024viscosity}.)  In particular, fix $\lambda^*:=4+\lfloor\tfrac{d}{2}\rfloor$, we can define a penalty function $ G(\mu)=G(\m; \mt):=\rho_{\lambda^*}(\mu,\mt)^2\ge0$.
		Then $G$ is $L$-differentiable with
		\begin{equation}
			\label{eq:G-Lions-derivative}
			\partial_\mu G(\mu)(x)=2\int_{\hR^d} \frac{-i\xi\, e^*(x,\xi)}{(1+|\xi|^2)^{\lambda^*}}\left(\mathfrak{F}(\mu)(\xi)-\mathfrak{F}(\mt)(\xi)\right)\d\xi,\qquad x\in \hR^d,
		\end{equation}
		and $G(\mu)=0$ if and only if $\mu=\mt$. That is, $G$ satisfies Definition~\ref{penalty}.
		%
		Moreover, from \eqref{eq:G-Lions-derivative} we can easily check, for $x\in \hR^d$ and $\mu\in \sP_2(\hR^d)$ it holds that
		\bea 
		\label{G-LD-bound}        
		|\partial_\mu G(\mu)(x)|\le 2\int_{\hR^d} \frac{2|\xi|}{(1+|\xi|^2)^{\lambda^*}}\d \xi< \infty,\quad    
		\|\nabla_x \partial_\mu G(\mu)(x)\|_{\infty}\le2\int_{\hR^d}\frac{2|\xi|^2}{(1+|\xi|^2)^{\lambda^*}}\d \xi< \infty,
		\eea
		which in particular implies that $G$ satisfies Assumption~\ref{assump1}-(2).    
		\qed
	\end{Example}
	
	
	
	\begin{Remark}
		{\rm We note that both the Gaussian MMD and Fourier--Wasserstein penalty functions are quadratic in $\m$, therefore convex under linear mixtures. {\it But they do not
				satisfy lifted convexity in general}. Indeed, let $\mu_{\mathrm{tar}}=\delta_0$. In the case of Gaussian MMD, we have
			$$
			G_{\mathrm{MMD}}(\delta_x;\delta_0)=2\Big(1-\exp\Big(-\frac{|x|^2}{2\eta^2}\Big)\Big).
			$$
			If the lift of $G_{\mathrm{MMD}}$ were convex, then its restriction to
			deterministic random variables would imply that
			$x\mapsto G_{\mathrm{MMD}}(\delta_x;\delta_0)$ is convex. Let $x=re$ for a fixed unit vector $e$, we have 
			$$
			\frac{\d^2}{\d r^2}
			G_{\mathrm{MMD}}(\delta_{re};\delta_0)
			= \frac{2}{\eta^2} \exp\Big(-\frac{r^2}{2\eta^2}\Big)\Big(1-\frac{r^2}{\eta^2}\Big)\le 0, \q |r|>\eta.
			$$
			Thus the $G_{\rm MMD}$  is not globally convex in the lifted sense.
			
			Similarly, for the Fourier--Wasserstein penalty, up to the fixed
			positive normalization constant in the definition of the Fourier
			transform,
			$$
			G_{\mathrm{F}}(\delta_x;\delta_0)
			=
			2
			\int_{\mathbb R^d}
			\frac{1-\cos(\xi\cdot x)}
			{(1+|\xi|^2)^{\lambda^\ast}}
			\,\d\xi.
			$$
			This is a bounded and nonconstant function of $x$. Since a finite
			convex function on all of $\mathbb R^d$ that is bounded above must be
			constant, the Fourier--Wasserstein penalty is also not globally
			convex in the lifted sense.
			\qed}
	\end{Remark}

	\paragraph{FBSDE formulation.}
	
	To end this section we shall derive the  McKean--Vlasov FBSDE {\color{black}corresponding to Problem 1 for each given $k>0$}, following the   well-known stochastic maximum principle for mean field systems (see, e.g.,  \cite{buckdahn2016stochastic} or \cite[Volume I, \S 6.2]{carmona2018probabilistic}). To this end, we consider the optimal control problem  \eqref{RSB_objective} and \eqref{eq:controlled-dynamics-follmer} , and define the corresponding Hamiltonian function $H$:
	\begin{equation}
		\label{eq:hamiltonian}
		H(t,x,y,\alpha) =\langle b(t,x)+\sigma(t)\alpha, y \rangle  + \frac{1}{2} |\alpha|^2,
	\end{equation}
	for $(t,x,y,\alpha)\in [0,T]\times \hR^d\times \hR^d\times \hR^d$.

	{\color{black} Now let $k>0$ be fixed. Applying the stochastic Pontryagin maximum principle} (see, \cite{buckdahn2016stochastic}
	and \cite[Theorem 6.14, 6.16]{carmona2015probabilistic}), and that the optimal control is of the form $\h\alpha(t,X_t,Y_t):=-\si(t)^TY_t$, $t\in[0,T]$, where $Y$ is the so-called \textit{adjoint process} and $X$ is the corresponding optimal trajectory, the closed-loop optimal system satisfies the following FBSDE:
	\begin{align}
		\label{mkvfbsde}
		\begin{cases}
			\d X_t = [b(t,X_t)+ \sigma(t)\h\alpha(t)] \d t + \sigma(t) \d W_t = [b(t,X_t)-\sigma(t)\sigma(t)^T Y_t]\d t + \sigma(t)\d W_t, \medskip\\
			\d Y_t = - (\nabla b(t,X_t))^TY_t \d t+ Z_t \d W_t,\medskip\\
			X_0 \sim \mu_{\rm ini}, \quad Y_T = k\partial_\mu G(\hP_{X_T})(X_T), \qq  t\in[0,T].
		\end{cases}
	\end{align}
	%
	Clearly, the system (\ref{mkvfbsde}) is a McKean--Vlasov type  FBSDE (or mean field FBSDE), whose well-posedness will be the main topic of the next section. 
	
	In what follows we denote $\theta_t:=(X_t,Y_t,Z_t)$, $t\in[0,T]$, to be triplet of progressively measurable processes  taking values in $\hR^d \times \hR^d\times \hR^{d\times d}$, such that $X, Y\in \hC([0,T];\hR^d)$, a.s. We denote $\Theta $ to be the collection of such processes $ \{\theta_t\}$, $t\in[0,T]$, and endow the norm  
	on the vector space $\Theta$: 
	\begin{equation}\label{def:norm_theta}
		\|\theta\|_\Theta = \hE\Big[ \sup_{0\leq t\leq T} \{|X_t|^2+ |Y_t|^2 \} + \int_{0}^T |Z_t|^2 \d t\Big]^{1/2}.
	\end{equation} 
	For notation simplicity, in what follows we denote 
	\bea
	\label{phimu}
	\phi_\mu(x):=\frac{\delta G}{\delta\mu} (\mu)(x),\qquad\forall\, \mu\in \sP_2(\hR^d),\,x\in \hR^d.
	\eea

	\section{Well-posedness of the optimal closed-loop FBSDE}

	In this section we prove the existence and uniqueness of the McKean-Vlasov type FBSDE \reff{mkvfbsde}. We note that due to the nonlinearity and the minimum regularity assumptions on the drift coefficient $b$, to the best of our knowledge, the well-posedness of such FBSDE is not covered by any existing literature.  Our proof relies on a {\it fixed-point}-type argument by first ``freeze'' the terminal distribution, and the existence and uniqueness of the minimizing path-law of Problem 2.

	We begin with the following lemma. 
	\begin{Lemma}[Existence and uniqueness of the optimal path law]
		\label{lem:path-law-existence-uniqueness}
		Under Assumptions~\ref{assump1}, the path law formulation of SCSBP \eqref{eq:path-law-objective} admits a unique minimizer $\hQ^{k,*}\in \sP(\Omega)$.
	\end{Lemma}
	The proof utilizes the lower semicontinuity of  $D(\,\cd\,\|\hP)$ and $G$. See details in Appendix~\ref{sec:appendix}.

	\subsection{Conditional Gibbs principle and self-consistency}
	
	We now give a detailed description of our ``freezing'' argument. To begin with, let us introduce the following auxiliary {\it frozen linearized functional} which will play an important role in the scheme.  Recall the linear derivative $\phi_\m(\cd)$ defined by \reff{phimu}.
	For fixed $\mu\in\sP_2(\mathbb R^d)$, consider a functional  $\cL_\m:\mathcal Q_{\mu_{\mathrm{ini}}}\mapsto \hR_+$ defined by 
	\bea
	\label{eq:frozen-linearized-functional}
	\mathcal L_\mu(\hQ) :=  D_{\mathrm{KL}}(\hQ\Vert \hP)
	+k\int_{\mathbb R^d}\phi_\mu(x)\hQ_T(\d x)=D_{\rm KL}(\hQ\|\hP)+k\hE^{\hQ}[\phi_\m(X_T)],  \qquad
	\hQ\in\mathcal Q_{\mu_{\mathrm{ini}}}, 
	\eea
	where $X$ is the canonical process and $\hQ_T=\hQ_{X_T}$. 
	Next, let $\hP\in\sQ_{\rm ini}$ be the path law of the reference diffusion (that is, under $\hP$, $X$ satisfies the SDE define by (\ref{P0})),  and denote the transition density of $X$ under $\hP$ by $p(t,y;s,x)$, $0\le s\le t\le T$, $x,y\in\hR^d$. Consider the corresponding transition operator:
	\bea
	\label{Pb}
	P_{s,t}[\f](x)=\int_{\hR^d} \f(y)p(t,y;s,x)dy=\hE^\hP[\f(X_t)|X_s=x], \q 0\le s\le t\le T, ~ \f\in\hC_b(\hR^d).
	\eea
	Now, define
	\begin{equation}
		\label{eq:conditional-normalization}
		Z_\mu(x)  :=  P_{0,T}\big[e^{-k\phi_\mu}\big](x)  =
		\mathbb E^\hP\Big[e^{-k\phi_\mu(X_T)}\Big|X_0=x\Big].
	\end{equation}
	Since $\phi_\mu$ is uniformly bounded, thanks to Assumption~\ref{assump1}, $Z_\mu$ is bounded above and bounded away from zero, and we can define a finite measure on $(\O, \cF)$ by 
	\bea
	\label{Qmu}
	\h\hQ^\mu(A)=\int_A  \frac{e^{-k\phi_\mu(X_T(\o))}}{Z_\mu(X_0(\o))}\hP(\d\o), \q A\in\cF.
	\eea
	In particular, if $A\in\si(X_0)$, that is, $A=X^{-1}_0(B)$ for some $B\in\sB(\hR^d)$, then by definition of conditional expectation we have
	\begin{equation}
		\h\hQ^\mu(A)=\int_{B}\frac{\hE^\hP[e^{-k\phi_\m(X_T)}|X_0=x]}{Z_\m(x)} \hP\circ X^{-1}_0(\d x)=\hP\{X_0\in B\} =\hP (A).
	\end{equation}
	This implies that $\h\hQ^\m(\O)=1$, whence $\h\hQ^\m \in\sP(\O)$ and is equivalent to $\hP$ with 
	$\h\hQ^\mu_{X_0}=\mu_{\mathrm{ini}}$. In other words, 
	$\h\hQ^\mu\in\mathcal Q_{\mu_{\mathrm{ini}}}$. 
	
	The following result characterize the nature of $\h\hQ^\m$.

	\begin{Lemma}[Conditional Gibbs variational principle]
		\label{lem:conditional-gibbs}
		For every $\mu\in\sP_2(\mathbb R^d)$, the law $\h\hQ^\mu$ is the unique minimizer of \eqref{eq:frozen-linearized-functional}. More precisely, for every $\h\hQ\in\mathcal Q_{\mu_{\mathrm{ini}}}$, it holds that
		\begin{equation}
			\mathcal L_\mu(\h\hQ)-\mathcal L_\mu(\h\hQ^\mu)
			=
			D_{\mathrm{KL}}(\h\hQ\Vert \h\hQ^\mu) \ge 0.
			\label{eq:conditional-gibbs-identity}
		\end{equation}
	\end{Lemma}
	
	\begin{proof}
		First note that the definition (\ref{Qmu}) amounts to saying that $\frac{d\h\hQ^\m}{d\hP}\big|_{\cF_T}=\frac{e^{-k\phi_\m(X_T)}}{Z_\m(X_0)}$. 
		Thus
		$\log\frac{d\h\hQ^\mu}{d\hP}=-k\phi_\mu(X_T)-\log Z_\mu(X_0)$. Next, note that  for every $\h\hQ\in\mathcal Q_{\mu_{\mathrm{ini}}}$, one has $\log\frac{d\h\hQ}{d\h\hQ^\m}= \log\frac{d\h\hQ}{d\hP}-\log\frac{d\h\hQ^\mu}{d\hP}$, and it then follows from  \reff{eq:frozen-linearized-functional} that 
		\begin{align}
			D_{\mathrm{KL}}(\h\hQ\Vert \h\hQ^\mu)
			&=
			D_{\mathrm{KL}}(\h\hQ\Vert \hP)
			+k\mathbb E^{\hat\hQ}[\phi_\mu(X_T)]
			+\mathbb E^{\hat\hQ}[\log Z_\mu(X_0)]=\cL_\m(\h\hQ)+\hE^{\hat\hQ}[\log Z_\mu(X_0)].
			\label{eq:entropy-change-reference}
		\end{align}
		In particular, this implies that 
		\bea
		\label{cLmu0}
		\mathcal L_\mu(\h\hQ^\mu)=-\mathbb E^{\h\hQ^\mu}[\log Z_\mu(X_0)].
		\eea
		But since $\h\hQ_{X_0}=\h\hQ^\mu_{X_0}=\mi$, their expectations of $\log Z_\mu(X_0)$ coincide. Plugging \reff{cLmu0} into \reff{eq:entropy-change-reference} yields \eqref{eq:conditional-gibbs-identity},
		which shows that $\h\hQ^\m$ is the minimizer of $\cL_\m(\cd)$. The equality condition for relative entropy gives uniqueness, completing the proof.
	\end{proof}
	
	Comparing the definition of the functionals $\cL_\m$ and $J^k$, it is conceivable that a ``fixed-point'' of the mapping $\m\mapsto \h\hQ^\m$ would give the minimizer of the original Problem 2 \reff{eq:path-law-objective}. The following result justifies this fact.  Recall the definition of Bregman remainder  $\mathcal B_G$ in Equation~\eqref{eq:linear-bregman-remainder}.

	\begin{Proposition}[Self-consistent Gibbs characterization]
		\label{prop:self-consistent-gibbs}
		A law $\hQ^*\in\mathcal Q_{\mu_{\mathrm{ini}}}$ is the minimizer of \eqref{eq:path-law-objective} if and only if
		\begin{equation}
			\hQ^*=\h\hQ^\mu,
			\qquad
			\mu=\hQ^*_T.
			\label{eq:self-consistent-gibbs-relation}
		\end{equation}
		Furthermore, if $\hQ^*$ satisfies \eqref{eq:self-consistent-gibbs-relation}, then for every $ \hQ\in\mathcal Q_{\mu_{\mathrm{ini}}}$, it holds that
		\begin{equation}
			J^k( \hQ)-J^k(\hQ^*)
			=  D_{\mathrm{KL}}( \hQ\Vert \hQ^*)  +k\mathcal B_G( \hQ_T\,|\,\mu)\ge 0.
			\label{eq:objective-gap-identity}
		\end{equation}
		Finally, such a minimizer $\hQ^*$ is unique in $\sQ_{\rm ini}$.
	\end{Proposition}
	
	\begin{proof}
		Suppose first that $\hQ^*$ satisfies \eqref{eq:self-consistent-gibbs-relation}. Applying Lemma~\ref{lem:conditional-gibbs} and using $\hQ^*_T=\mu$ we have, for any $\hQ\in\sQ_{\rm ini}$, 
		\begin{align*}
			J^k( \hQ)-J^k(\hQ^*)
			&=
			\mathcal L_\mu(\hQ)-\mathcal L_\mu(\hQ^*) +k\Big[G( \hQ_T)-G(\mu) -\int_{\mathbb R^d}\phi_\mu(x)( \hQ_T-\mu)(\d x)\Big]\\
			&= D_{\mathrm{KL}}( \hQ\Vert \hQ^*)
			+k\mathcal B_G( \hQ_T\,|\,\mu).
		\end{align*}
		Since $D_{\mathrm{KL}}( \hQ\Vert \hQ^*)\ge0$, and $\mathcal B_G( \hQ_T\,|\,\mu)\ge 0$, thanks to Assumption~\ref{assump1}
		and Remark \ref{rem:linear-versus-lifted-convexity}, we obtain  \eqref{eq:objective-gap-identity}, which in turn shows that
		$\hQ^*$ is a minimizer of $J^k$ on $\sQ_{\rm ini}$.
		
		Conversely, let $\hQ^*$ minimize $J^k$ and define $\mu:=\hQ^*_T$. Let $\h\hQ:=\h\hQ^\mu$, the unique minimizer of $\mathcal L_\mu$. We claim that $\hQ^*$ must satisfy \reff{eq:self-consistent-gibbs-relation}. Indeed, suppose $\hQ^*\neq \h\hQ^\mu$. Then by the uniqueness of the minimizer we must have 
		$\mathcal L_\mu(\h\hQ)<\mathcal L_\mu(\hQ^*)$, we shall derive a contradiction.

		For $\varepsilon\in(0,1)$, define
		$
		\hQ^\varepsilon:=(1-\varepsilon)\hQ^*+\varepsilon \h\hQ$. 
		By the convexity of relative entropy we have 
		$$
		D_{\mathrm{KL}}(\hQ^\varepsilon\Vert \hP)
		\leq
		(1-\varepsilon)D_{\mathrm{KL}}(\hQ\Vert \hP)
		+\varepsilon D_{\mathrm{KL}}(\h\hQ\Vert \hP).
		$$
		Since $(\hQ^\varepsilon)_T=(1-\varepsilon)\mu+\varepsilon \h\hQ_T$,   the definition of the linear functional derivative then yields
		$$ G(\hQ^\varepsilon_T) = G(\mu) + \varepsilon\int_{\mathbb R^d}\phi_\mu(x)(\h\hQ_T-\mu)(dx)
		+o(\varepsilon).
		$$
		Consequently, it is readily seen that 
		$$
		J^k(\hQ^\varepsilon)-J^k(\hQ^*)
		\leq
		\varepsilon\bigl[\mathcal L_\mu(\h\hQ)-\mathcal L_\mu(\hQ^*)\bigr]
		+o(\varepsilon)<0,
		$$
		for $\varepsilon>0$ sufficiently small. This contradicts the optimality of $\hQ^*$, proving \eqref{eq:self-consistent-gibbs-relation}.
		
		It remains to prove the uniqueness. 
		Suppose that there are two minimizers, 
		$\hQ^i:=\hQ^{\mu_i}$ with $\mu^i= \hQ^i_T$, $i=1,2$. Since  $D_{\rm KL}(\hQ^i\| \hP)<\infty$ and $\frac{d\hO^i}{d\hP}$ are bounded both above and   away from zero, we have $D_{\rm KL}(\hQ^i\|\hQ^j)<\infty$, $i,j =1,2$. In particular, $\hQ^1$ and $\hQ^2$ are equivalent.
		
		Since both $\hQ^1$ and $\hQ^2$ are minimizers,  $J^k(\hQ^1)=J^k(\hQ^2)=V^k$. Thus by \reff{eq:objective-gap-identity}, we have
		\bea
		\label{eq:cross-objective-gap}
		0=J^k(\hQ^2)-J^k(\hQ^1)  =
		D_{\mathrm{KL}}(\hQ^2\|\hQ^1)  +k\mathcal B_G(\mu^2\,|\,\mu^1)\ge 0.
		\eea
		Again, note that both terms on the right-hand side above is nonnegative, we conclude that 
		both terms must be zero. Switching the position of $\hQ^1$ and $\hQ^2$ we obtain
		$$D_{\mathrm{KL}}(\hQ^i\|\hQ^j)=\mathcal B_G(\mu^i\,|\,\mu^j)=0, \qq i,j \in\{1,2\}.
		$$
		Therefore  $\hQ^1=\hQ^2$ and $\mu^1=\mu^2$. The proof is now complete. 
	\end{proof}

	\subsection{The Decoupling Fields}
	
	In this subsection we construct an important device for solving FBSDE \reff{mkvfbsde}, that is, the so-called {\it decoupling field}. 
	To begin with,  fix $\m\in\sP_1(\hR^d)$, and recall the normalizer $Z_\m(\cd)$ 
	defined by \reff{eq:conditional-normalization}. Consider its ``dynamic'' version:
	\bea
	\label{hmu}
	h^\mu(t,x) := \mathbb E^\hP\Big[e^{-k\phi_\mu(X_T)}\Big| X_t=x\Big], \qq (t,x)\in[0,T]\times\hR^d.
	\eea
	In what follows we denote $q_\mu(x):=e^{-k\phi_\mu(x)}$ for simplicity. Then by Assumption~\ref{assump1}(2)(iv), we have
	\bea
	\label{qmu}
	e^{-kL_G}\leq q_\mu(x)\leq e^{kL_G}, \q  \sup_{\mu\in\sP_2(\mathbb R^d)} \sum_{j=1}^3\|\nabla_x^jq_\mu\|_\infty \leq L_G.
	\eea
	Now we consider the reference diffusion on canonical space $(\O, \cF, \hP^0)$, where $\hP^0$ is Wiener measure:
	\begin{equation}
		dX_s^{t,x} =  b(s,X_s^{t,x})\,ds+\sigma(s)\,dW_s, \qq s\in[t,T],  \quad  X_t^{t,x}=x.
		\label{eq:reference-diffusion-starting-tx}
	\end{equation}
	Clearly, under Assumption~\ref{assump1}, SDE \reff{eq:reference-diffusion-starting-tx} has a unique strong solution, and we denote its law under $\hP^0$ by $\hP$, as before. In light of the Feynman-Kac formula, we know that $h^\m$ satisfies the following PDE:
	\bea
	\label{eq:frozen-kolmogorov}
	\left\{\ba{lll}
	\dis  \partial_th^\mu(t,x)  +b(t,x)\cdot\nabla_xh^\mu(t,x) +\frac12 a(t):\nabla_x^2h^\mu(t,x)  =0, \q (t,x)\in [0,T]\times \hR^d; \\
	h^\mu(T,x)=q_\mu(x),
	\ea\right.
	\eea
	where $a(t)=\si(t)\si^\top(t)$ and $A:B:=\operatorname{tr}(AB)$, for $A,B\in\hR^{d\times d}$. Furthermore, given the regularity of $q_\m$ (\ref{qmu}), we know that $h^\m$ is the unique classical solution, such that $h^\m$, $\nabla h^\m$, and $\nabla^2 h^\m$ are all uniformly bounded, uniformly in $\m\in \sP_2(\hR^d)$. In particular, the $\hC^3$-regularity of $q_\m$ renders that $h^\m\in\hC^{1,3}([0,T]\times \hR^d)$ (see, e.g., \cite{Lady1968} or \cite{Lieberman1996}). 
	
	Let us now consider the following logarithmic   (or Hopf-Cole) transform:
	\begin{equation}
		\label{eq:frozen-log-transform}
		v^\mu(t,x) :=  -\log h^\mu(t,x), \qquad u^\mu(t,x) := \nabla_xv^\mu(t,x).
	\end{equation}
	We note that  the maximum principle, together with the property of $q_\m$ (\ref{qmu}) guarantees that $h^\mu(t,x)\geq e^{-kL_G}>0$, 
	$(t,x)\in [0,T]\times \hR^d$ (cf. e.g., \cite{Friedman1964}). Thus $v^\m$ is well-defined and $v^\m\in\hC^{1,3}([0,T]\times \hR^d)$ as well, and consequently, $v^\m\in\hC^{1,2}([0,T]\times\hR^d)$. Moreover, it is easy to check that $v^\m$ satisfies the PDE:
	\bea
	\label{eq:frozen-hjb}
	\left\{\ba{lll}
	\dis \partial_tv^\mu(t,x)
	+b(t,x)\cdot\nabla_xv^\mu(t,x)
	+\frac12 a(t):\nabla_x^2v^\mu(t,x) -\frac12\left|\sigma(t)^\top\nabla_xv^\mu(t,x)\right|^2
	=0,\\
	v^\mu(T,x)=k\phi_\mu(x);
	\ea\right.
	\eea
	whereas $u^\m$ satisfies the PDE: 
	\bea
	\label{eq:differentiated-frozen-hjb}
	\left\{\ba{lll}
	\dis \partial_tu^\mu+(\nabla_xb)^\top u^\mu
	+\nabla_xu^\mu \cd (b-au^\mu)
	+\frac12 a(t): \nabla_x^2u^\mu  =0,\\
	u^\mu(T,x)=k\nabla_x\phi_\mu(x),
	\ea\right.
	\eea
	where all coefficients in the first line are evaluated at $(t,x)$. Notice that $a(t):\nabla_x^2u^\mu$ is taken componentwise, so that $\bigl(a(t):\nabla_x^2u^\mu\bigr)_i
	= a(t):\nabla_x^2u_i^\mu
	= \sum_{j,\ell=1}^{d}a_{j\ell}(t)\,\partial_{x_jx_\ell}u_i^\mu$ for  $i=1,\ldots,d$.

	We summarize all the facts above in the following lemma for ready reference. 
	\begin{Lemma}
		\label{lem:frozen-log-transform-regularity}
		Assume Assumptions~\ref{assump1} is in force. Then, for every $\mu\in\sP_2(\mathbb R^d)$,
		$v^\mu\in \hC^{1,3}([0,T]\times\mathbb R^d)$ is the classical solution of the PDE \reff{eq:frozen-hjb}, and
		$u^\mu=\nabla_xv^\mu\in \hC^{1,2}([0,T]\times\mathbb R^d;\mathbb R^d)$ is the classical solution to the PDE \reff{eq:differentiated-frozen-hjb}.
		
		Moreover, there exists a finite constant $C_{T,k}$, independent of $\mu$, such that
		\begin{equation}
			\sup_{\mu\in\sP_2(\mathbb R^d)}
			\sup_{(t,x)\in[0,T]\times\mathbb R^d}
			\left(
			|u^\mu(t,x)|+|\nabla_xu^\mu(t,x)|
			\right)
			\leq C_{T,k}.
			\label{eq:uniform-frozen-field-bound}
		\end{equation}
	\end{Lemma}

	\subsection{Well-posedned of the Closed-loop FBSDE (\ref{mkvfbsde})}
	
	We now turn to the well-posedness of the McKean-Vlasov FBSDE \reff{mkvfbsde}. To begin with, we   consider the following ``frozen'' version of it for a given $\m\in\sP_2(\hR^d)$, which is of a standard Markovian type on any filtered probability space: 
	\begin{equation}
		\label{eq:frozen-fbsde}
		\left\{\ba{lll} 
		\d X_t  =  \left(b(t,X_t)-a(t)Y_t\right)\d t  +\sigma(t)\,\d W_t, \q & X_0 =\xi_0, \ms \\
		\d Y_t  = -\left(\nabla_xb(t,X_t)\right)^\top Y_t\,\d t  +Z_t\,\d W_t,  &Y_T =k\nabla_x\phi_\mu(X_T), 
		\ea \right.
	\end{equation}
	where $\xi_0\sim \mi$ and $a(t)=\si(t)\si^\top(t)$. In particular, we shall consider the canonical space $(\O, \cF, \hP)$, where $\hP$ is the path-law of the reference diffusion \reff{P0}. The following result is  more or less standard by the Four-Step Scheme for FBSDEs, given the decoupling fields $h^\m$ and $u^\m$  (cf. e.g., \cite{MWZZ2015,MaYong1999}). 
	\begin{Lemma}[Well-posedness of the frozen FBSDE]
		\label{lem:frozen-fbsde}
		Assume that Assumption~\ref{assump1} is in force. Fix $\mu\in\sP_2(\mathbb R^d)$, and $\xi_0\in\hL^2(\O)$ that is independent of the $\hP$-Brownian motion $W$ with 
		$ 
		\hP_{\xi_0}=\mu_{\mathrm{ini}},
		$ 
		the FBSDE \reff{eq:frozen-fbsde} has a unique strong solution  $\Theta=(X, Y, Z)$ on $(\O, \cF, \hP)$.
		Moreover, the solution is given by
		\begin{equation}
			Y_t=u^\mu(t,X_t),
			\qquad
			Z_t=\nabla_xu^\mu(t,X_t)\sigma(t),
			\label{eq:frozen-decoupling-representation}
		\end{equation}
		and the path law of $X$ is $\h\hQ^\mu$.
	\end{Lemma}
	
	\begin{proof}
		We first follow the Four Step Scheme to show that the FBSDE (\ref{eq:frozen-fbsde}) has the {\it nodal solution} $\tilde\Th=(\tilde X, \tilde Y, \tilde Z)$, where  $\tilde Y$ and $\tilde Z$ are functions of $\tilde X$ defined by \eqref{eq:frozen-decoupling-representation}, and $\tilde X$ satisfies the 
		decoupled SDE:
		\begin{equation}
			\d\tilde{X}_t = \left(b(t,\tilde X_t)-a(t)u^\mu(t,\tilde X_t)\right)\d t  +\sigma(t)\,
			\d W_t, \q \tilde X_0=\xi_0.
			\label{eq:frozen-forward-sde}
		\end{equation}
		
		Indeed, by \eqref{eq:uniform-frozen-field-bound}, the function $ u^\mu(t,\cd)$ is uniform Lipschitz continuous,  uniformly in $t\in[0,T]$. Therefore  \eqref{eq:frozen-forward-sde} has a unique strong solution $\tilde X$ on $(\O, \cF, \hP)$, thanks to Assumption \ref{assump1} on the coefficients $b$ and $\si$. Next, defining $(\tilde Y, \tilde Z)$ by \reff{eq:frozen-decoupling-representation}, applying It\^o's formula to $u^\mu(t,\tilde X_t)$, and noting that $u^\m$ satisfies the PDE
		\eqref{eq:differentiated-frozen-hjb}, one can easily check that 
		$$ \d\tilde Y_t = -\left(\nabla_xb(t,\tilde X_t)\right)^\top \tilde Y_t\,\d t +\tilde Z_t\,\d W_t, \q 
		\tilde Y_T =k\nabla_x\phi_\mu(\tilde X_T).
		$$
		That is, $(\tilde X, \tilde Y, \tilde Z)$ satisfies the FBSDE \reff{eq:frozen-fbsde}. Furthermore, given the regularity of the coefficients $b$, $\si$ and the decoupling field $u^\m$, it can be shown (cf. \cite[Chapter 4, Theorem 1.1]{MaYong1999} that every adapted solution of \reff{eq:frozen-fbsde} must coincide with the nodal solution in the pathwise sense. This proves the strong well-posedness of \reff{eq:frozen-fbsde} on $(\O, \cF, \hP)$ and the decoupling representation \reff{eq:frozen-decoupling-representation}. 
		%
		
		It remains to show that 
		path law of $\tilde X$ is $\h\hQ^\mu$. To this end, note that under $\hP$, the canonical process $X$ satisfies the SDE (\ref{P0}). Define
		\begin{equation}
			M_t^\mu  :=
			\frac{h^\mu(t,X_t)}{h^\mu(0,X_0)}, \qq t\in[0,T].
			\label{eq:gibbs-density-martingale}
		\end{equation}
		Applying It\^o's formula to $h^\mu(t,X_t)$ and using \eqref{eq:frozen-kolmogorov} we see that $M^\m$ satisfies the linear SDE:
		\begin{equation}
			\d M_t^\mu
			= -M_t^\mu \cd  (\sigma(t)^\top u^\mu(t,X_t),\d W_t ), \q t\in[0,T].
			\label{eq:gibbs-stochastic-exponential}
		\end{equation}
		Moreover, recall the probability measure $\h\hQ^\m$  (\ref{Qmu}) and the normalizer $Z_\m$ \reff{eq:conditional-normalization}, we see that  
		$ 
		M_T^\mu= \frac{e^{-k\phi_\mu(X_T)}}{Z_\mu(X_0)}=\frac{\d\h\hQ^\mu}{\d\hP}.
		$ 
		By Girsanov's theorem, we conclude that under $\h\hQ^\mu$, the (canonical) process $X$ satisfies the SDE
		\begin{equation}
			\d X_t  =  \left(b(t,X_t)-a(t)u^\mu(t,X_t)\right)\d t
			+\sigma(t)\,\d W_t^\mu, \q t\in[0,T],
			\label{eq:gibbs-controlled-dynamics}
		\end{equation}
		where $W^\mu_t:=W_t+\int_0^t\si(t)^\top u^\m(s, X_s)\d s$ is a $\h\hQ^\mu$-Brownian motion. Clearly, the SDE (\ref{eq:gibbs-controlled-dynamics}) and (\ref{eq:frozen-forward-sde}) have the same form, thus the weak uniqueness implies that $\tilde X$ and $X$ have the same path law, which is $\h\hQ^\m$ since  $X$ is the canonical process. The proof is now complete. 
	\end{proof}
	
	%
	
	We are now ready to prove the well-posedness of the McKean--Vlasov FBSDE \reff{mkvfbsde}. Let the time duration $T>0$ be given, and still consider the canonical space $(\O, \cF, \hP)$, where $\hP$ is the path law of the reference diffusion \reff{P0}. . 
	
	\begin{Theorem}[Global well-posedness]
		\label{thm:global-fbsde-wellposedness}
		Assume Assumption~\ref{assump1} is in force. Then,
		
		\ms
		(1) for  any $k>0$, the McKean--Vlasov FBSDE \eqref{mkvfbsde} has  a solution $\Th^k=(X^k ,Y^k ,Z^k )$ defined on $(\O, \cF, \hP)$, and
		the path law of $X^k$ under $\hP$ is the unique minimizer $\hQ^{k,*}$ of \eqref{eq:path-law-objective}; 
		
		\ms
		
		(2) the solution of \reff{mkvfbsde} is unique in law; and 
		
		\ms
		(3)   the optimal control  of Problem 1 is given by 
		$ \widehat\alpha_t^k:=  -\sigma(t)^\top Y_t^k$.
	\end{Theorem}
	
	\begin{proof}
		(1) First recall from Lemma~\ref{lem:path-law-existence-uniqueness} that the path-law Problem 2 has a unique minimizer $\hQ^{k,*}$. Let $\mu_T^{k,*}:=\hQ^{k,*}_T$. Then by Proposition~\ref{prop:self-consistent-gibbs} we must have
		$ \hQ^{k,*}=\widehat \hQ^{\mu_T^{k,*}}$.
		
		Next, applying Lemma~\ref{lem:frozen-fbsde} with $\mu=\mu_T^{k,*}$  we obtain a unique strong solution of the frozen FBSDE \reff{eq:frozen-fbsde} on $(\O, \cF, \hP)$, such that 
		the path law of $X^k$ under $\hP$ is $\hP\circ (X^k)^{-1}=\h\hQ^{\mu_T^{k,*}}=\hQ^{k,*}$. In particular, we have 
		$\hP_{X_T^k}=\hQ^{k,*}_T=\mu_T^{k,*}$, 
		and the   terminal condition of $Y^k$ becomes
		$$ Y_T^k = k\nabla_x\phi_{\mu_T^{k,*}}(X_T^k) = k\partial_\mu G(\hP_{X_T^k})(X_T^k).
		$$
		In other words, $(X^k,Y^k,Z^k)$ actually solves the McKean-Vlasov FBSDE \eqref{mkvfbsde}, proving (1).
		
		(2) We shall argue directly from the linear convexity. Let $(X^i,Y^i,Z^i)$, $i=1,2$, be two square-integrable solutions, possibly defined on different stochastic bases $(\O^i, \cF^i, \hP^i, \hF^i, \xi_0^i, W^i)$, such that $\xi_0^i\sim \mi$, and $W^i$ is $\hP^i$ Brownian motion, $i=1,2$. Set
		$$
		\hQ^i:=\hP^i\circ (X^i)^{-1},
		\qquad
		\mu_i:=\hP^i_{X_T^i}.
		$$
		Since each $(X^i,Y^i,Z^i)$ must be the nodal solution of the  frozen FBSDE \eqref{eq:frozen-fbsde}, $X^i$ must satisfy the decoupled SDE 
		\reff{eq:frozen-forward-sde} (driven by $W^i$). Uniqueness in law for the SDE and Lemma~\ref{lem:frozen-fbsde} therefore imply that 
		$\hQ^i=\h\hQ^{\mu_i}$, whence $\hQ^i_T=\mu_i$, $i=1,2$. Thus both pairs $(\hQ^i,\mu^i)$ are self-consistent. Then it follows from Proposition
		\ref{prop:self-consistent-gibbs} that $\hQ^i$'s are the minimizers of Problem \eqref{eq:path-law-objective}. But by Lemma \ref{lem:path-law-existence-uniqueness}, such minimizer is unique. Thus we must have
		$$\hQ^1=\hQ^2=\hQ^*, 
		\qquad\text{and}\qquad
		\mu^1=\mu^2=:\mu^*.
		$$
		In other words, both $(X^i,Y^i,Z^i)$'s are the solutions of the same frozen FBSDE with parameter $\mu^*$. Lemma~\ref{lem:frozen-fbsde} gives uniqueness in law of the full triples.
		
		
		(3) Finally, every solution of \eqref{mkvfbsde} has a self-consistent Gibbs path law. Proposition~\ref{prop:self-consistent-gibbs} therefore shows that its path law minimizes \eqref{eq:path-law-objective}. Under the minimizing law, the forward equation has drift modification
		$$
		\sigma(t)\widehat\alpha_t^k=-a(t)Y_t^k,
		$$
		proving (3). 
		Uniqueness of the path law and uniqueness of the finite-entropy drift imply uniqueness of the optimal canonical drift $\dd t\otimes \dd\hQ^{k,*}$-almost everywhere.
	\end{proof}

			\section{Further Approximations}
			An important feature of the soft-constrained SBP is its convergence as $k\to \infty$. \cite{ma2025schr} proved that the convergence rate is ``linear'': for any $R>0$, there exists $C_R>0$ such that
			\beaa
			\int_0^T|\h\a^k(t,x)-\h \a(t,x)|dt \le \frac{C_R}{k}, \qq (t, x)\in[0,T]\times B_R, 
			\eeaa
			where $\h \a^k$ and $\h\a$ are the optimal controls of the SCSBP and the original SBP, respectively, and $B_R\subset\hR^d$ is the closed ball with radius $R$; see \cite[Theorem 6.4]{ma2025schr}. This result, however, is theoretical and does not address implementation. In particular, the SCSBP control $\h\a^k$ remains in the form of a ``score function'' and is no easier to simulate than the original control $\h \a$. In this section, we focus on data-driven approximations that address these implementation challenges.
			
			There are two main challenges arise. First, the distribution $\mt$ is both the object to be learned through the SCSBP and an input to the penalty function $G(\cd)=G(\cd;\mt)$ (i.e., $G(\mt)=0$) in Definition \ref{penalty}, and hence to Problem \ref{softSBP}. This circularity creates a fundamental implementation difficulty. Second, solving the SCSBP requires simulating the law $\hP_{X^\a_T}$, since this law is generated by an optimal closed-loop McKean--Vlasov FBSDE, its simulation is generally challenging.
			
			To address these issues, the next two sections introduce two auxiliary approximation problems. We first approximate $\mt$ by an empirical measure constructed from data sampled under $\mt$. We then approximate $\hP_{X_T^\alpha}$ using interacting particles to handle the mean-field nature of the FBSDE.

			\subsection{Data sample approximation}
			\label{sec:data_approximation}
			We begin by approximating the target distribution $\mt$ by an empirical distribution $\widehat{\m}^M$ defined by 
			$\widehat\mu^M:= \frac1M\sum_{j=1}^M \delta_{\zeta^j}$,
			where $\zeta^1,\ldots,\zeta^M$ are i.i.d. copies from $\mu_{\rm tar}$, independent of $(X_0,W)$. 
			
			Noting that the empirical measure $\widehat{\m}^M$ is by definition random, in light of Definition \reff{penalty} we consider the sample-averaged empirical approximation of $G(\mu):= G(\mu;\mu_{\rm tar})$ defined  by
			\bea
			\label{eq:G^M}
			G^M(\mu) :=  \hE\big[ G(\mu;\widehat\mu^M) \big],  \qquad \mu\in \sP_2(\mathbb R^d).
			\eea
			Clearly, $G^M(\cd)$ is a deterministic functional on $\sP_2(\mathbb R^d)$.
			We now replace the penalty function $G(\cd)$ in Problem \ref{softSBP} by $G^M(\cd)$, and consider the corresponding approximated problem.  Based on our previous theoretical discussion, it is reasonable to expect that the optimal control will be in the feedback form of the states. Therefore, in what follows we shall focus on  the closed-loop admissible controls in the form $\a_t=\a(t,X_t)$, $t\in[0,T]$, where $\a:[0,T]\times\mathbb{R}^d \to \hR^d$ is the  deterministic measurable functions. We shall  
			denote $\mathbb{A}$ to be the set of all such functions. Consequently, we now consider the first approximating version of Problem \ref{softSBP}:
			\begin{Problem}
				\label{prob3}
				Given $k, M>0$ and penalty function $G$,  find $V^{k,M}: =\inf_{\alpha \in \mathbb{A}} J^{k,M}(\alpha) $, where
				\begin{align}
					\label{eq:relaxed_SB-objective}
					\begin{split}
						J^{k,M}(\alpha)&:= 
						\hE \left[ \int_0^T \frac{1}{2}|\alpha(s,X^\a_s)|^2 \d s  + k G^M(\hP_{X^\alpha_T})\right], 
					\end{split}
				\end{align}
				and 
				\begin{equation}\label{eq:relaxed_SB-dynamics1}
					\d X^\alpha_t = \big(b(t,X^\alpha
					_t) + \sigma(t) \alpha(t,X^\a_t)\big)\d t + \sigma(t)\d W_t, \quad \hP_{X_0^\alpha} = \mu_{\rm ini}.	
				\end{equation}
				Here in  \eqref{eq:relaxed_SB-objective}, the expectation 
				is taken respect to any $\hP\in\sP(\O)$ under which  $W=(W_t)_{t \ge 0}$ is a Brownian motion  and the initial distribution $\hP_{X_0}=\mu_{\rm ini}$.
				\qed
			\end{Problem}

			Similar to Problem \ref{softSBP}, we see that the counterpart of FBSDE corresponding to Problem \ref{prob3} should take the following form:
			\begin{align}
				\label{eq:mkvfbsde_relaxed_km}
				\begin{cases}
					\d X^{k,M}_t = \big((b(t,X^{k,M}_t)- \sigma(t)\si(t)^\top Y^{k,M}_t\big) \d t + \sigma(t) \d W_t, \medskip\\
					\d Y^{k,M}_t = - \nabla b(t,X^{k,M}_t)Y^{k,M}_t \d t+ Z^{k,M}_t \d W_t,\medskip\\
					X^{k,M}_0 \sim \mu_{\rm ini}, \quad Y^{k,M}_T = k\partial_\mu G^M(\hP_{X^{k,M}_T})(X^{k,M}_T),
				\end{cases}
				\quad t\in[0,T].
			\end{align}
			
			Later we show  that  Problem~\ref{prob3} provides tractable approximation to Problem~\ref{softSBP}. To facilitate our discussion, let us make the following assumptions about the (approximate) penalty function $G^M$. 
			
			\begin{Assumption}
				\label{assum:data_sample}
				We assume that for any $M>0$, the functional $G^M(\cd)$  defined by \eqref{eq:G^M} satisfies the following properties:
				
				\bs
				(i)  $    \partial_\mu G^M(\mu)(x)=\hE\left[\partial_\mu G(\mu;\widehat\mu^M)(x)\right]$, for every $(\mu,x)\in \sP_2(\mathbb R^d)\times\mathbb R^d$.
				
				\bs
				
				(ii) There exists a constant $L_0>0$, independent of $M$, such that 
				$$
				\sup_{\mu\in\sP_2(\mathbb R^d)}\|\partial_\mu G(\mu)(\cdot)\|_\infty+\sup_{M\ge1}\sup_{\mu\in\sP_2(\mathbb R^d)}
				\|\partial_\mu G^M(\mu)(\cdot)\|_\infty  \le L_0,
				$$
				and  for all $X,X'\in L^2(\Omega;\mathbb R^d)$ and all $M\ge1$, it holds that
				$$
				\left\{\ba{lll}
				\mathbb E\left[\left|\partial_\mu G(\hP_{X'})(X')-\partial_\mu G(\hP_{X})(X)\right|^2\right]\le L_0\mathbb E[|X'-X|^2], \ms\\
				\mathbb E\left[\left|\partial_\mu G^M(\hP_{X'})(X')-\partial_\mu G^M(\hP_{X})(X)\right|^2\right]\le L_0\mathbb E[|X'-X|^2],
				\ea\right.
				$$

				\ms
				(iii)  For every $R>0$, define, for each $M>0$, 
				\bea
				\label{eoMR}
				\left\{\ba{lll}
				\dis    \varepsilon_M(R):= \sup_{\mu\in\sP_2(\mathbb R^d):\,\int |x|^2\mu(dx)\le R} \left|G^M(\mu)-G(\mu)\right|, \ms\\
				\dis   \omega_M(R):=  \sup_{X\in L^2(\Omega;\mathbb R^d)\atop \mathbb E|X|^2\le R} \hE\left[\left|\partial_\mu G^M(\hP_X)(X)-
				\partial_\mu G(\hP_X))(X) \right|^2\right].
				\ea\right.
				\eea
				Then, for every $R>0$, it holds that 
				$\varepsilon_M(R)\rightarrow 0$,  $ \omega_M(R)\rightarrow 0$,    as $  M\to\infty$.
				\qed
			\end{Assumption}
			
			These assumptions are key to establishing the convergence rate of the optimal policy; see Proposition \ref{prop:sample_approx}. They are relatively mild and are satisfied by many commonly used penalty functions. To illustrate this, we verify them below for several standard examples, including those introduced in Examples \ref{ex:mmd-gaussian} and \ref{ex: fourier-wasserstein}.
			\begin{Example}
				\label{rem:GM_rates}
				(1) We first consider the ``Fourier-Wasserstein'' penalty function in Example \ref{ex: fourier-wasserstein}:
				$$ G(\mu;\nu)=\rho_{\lambda^\ast}^2(\mu,\nu)= \int_{\mathbb R^d}(1+|\xi|^2)^{-\lambda^\ast}\left|\mathfrak F(\mu)(\xi)- \mathfrak F(\nu)(\xi)\right|^2 \d\xi,
				$$
				where $\mathfrak{F}(\mu)(\xi):=\int_{\hR^d} (2\pi)^{-\frac{d}{2}}e^{i\xi\cdot x}\mu(d x)$, $\xi\in \hR^d$, is the Fourier transform of 
				$\mu\in \sP_2(\hR^d)$. Then, it is not hard to check that 
				$$    G^M(\mu)=\hE\Big[ \int_{\mathbb R^d}(1+|\xi|^2)^{-\lambda^\ast}\Big|\mathfrak F(\mu)(\xi)-\mathfrak F(\widehat\mu^M)(\xi)\Big|^2 \d\xi =G(\mu)+r_M,
				$$
				where $r_M:=\int_{\mathbb R^d}(1+|\xi|^2)^{-\lambda^\ast} \mathbb E \big[\big| \mathfrak F(\widehat\mu^M)(\xi)-\mathfrak F(\mu_{\rm tar})(\xi)\big|^2\big]\d\xi$. We remark that in the above calculation the cross term vanishes because
				$ \mathbb E[\mathfrak F(\widehat\mu^M)(\xi)]=\mathfrak F(\mu_{\rm tar})(\xi)$.

				To verify Assumption \ref{assum:data_sample}-(i) and (ii), set
				$c_d=(2\pi)^{-d/2}$ and $w(\xi)=(1+|\xi|^2)^{-\lambda^\star}$.
				Since $\lambda^\star>d/2+1$, differentiation under the integral gives
				$$\partial_\mu G(\mu;\nu)(x)
				=
				2c_d\operatorname{Re}\int_{\mathbb R^d}
				i\xi e^{i\xi\cdot x}
				\overline{\mathcal F(\mu)(\xi)-\mathcal F(\nu)(\xi)}
				w(\xi)\,\d\xi.$$
				
				The identity $G^M=G+r_M$ and the unbiasedness of
				$\mathcal F(\widehat\mu^M)$ therefore imply
				\[
				\partial_\mu G^M(\mu)(x)
				=\partial_\mu G(\mu;\mu_{\rm tar})(x)
				=\mathbb E\big[\partial_\mu G(\mu;\widehat\mu^M)(x)\big],
				\]
				proving (i).
				
				Using $|\mathcal F(\mu)|\le c_d$ and
				$|\mathcal F(\mu)(\xi)-\mathcal F(\mu')(\xi)|
				\le c_d|\xi|W_2(\mu,\mu')$, the derivative formula yields
				\[
				\begin{aligned}
					|\partial_\mu G(\mu;\nu)(x)|\le C, \quad
					|\partial_\mu G(\mu;\nu)(x)-\partial_\mu G(\mu';\nu)(x')|
					\le C\bigl(|x-x'|+W_2(\mu,\mu')\bigr),
				\end{aligned}
				\]
				where $C$ depends only on $d$ and $\lambda^\star$.
				Consequently, $W_2(\mathcal L(X),\mathcal L(X'))^2
				\le\mathbb E|X-X'|^2$ implies
				\[
				\mathbb E\left|
				\partial_\mu G(\mathcal L(X))(X)
				-\partial_\mu G(\mathcal L(X'))(X')
				\right|^2
				\le 4C^2\mathbb E|X-X'|^2.
				\]
				Since $\partial_\mu G^M=\partial_\mu G$, the same bounds hold for
				$G^M$, uniformly in $M$. Thus (ii) holds with
				$L_0=\max\{2C,4C^2\}$.

				We now check Assumption \ref{assum:data_sample}-(iii). Indeed, note that
				$\partial_\mu G^M(\mu)(x) = \partial_\mu G(\mu;\mu_{\rm tar})(x)$, we have $\omega_M(R)\equiv0$. Moreover, since 
				$
				r_M \le \frac{C}{M} \int_{\mathbb R^d}(1+|\xi|^2)^{-\lambda^\ast} \d\xi \le\frac{C}{M}, 
				$
				we conclude that $ \varepsilon_M(R)\le \frac{C}{M}$. In other words, for the averaged Fourier--Wasserstein penalty, the empirical target sampling changes the value by an additive $O(M^{-1})$ bias but does not change the optimal feedback control.
				
				\ms
				
				(2) We now look at the ``Gaussian MMD'' penalty function in Example \ref{ex:mmd-gaussian}:
				$$    G(\mu;\nu):=
				\iint k_\eta(x,x')\,\mu(\d x)\mu(\d x') -2\iint k_\eta(x,y)\,\mu(\d x)\nu(\d y)     +
				\iint k_\eta(y,y')\,\nu(\d y)\nu(\d y').
				$$
				A direct computation shows that a similar cancellation also occurs, and we have 
				$$
				G^M(\mu)= G(\mu;\mu_{\rm tar}) + \frac1M \left[ \int k_\eta(y,y)\,\mu_{\rm tar}(\d y) - \iint k_\eta(y,y')\,
				\mu_{\rm tar}(\d y)\mu_{\rm tar}(\d y')\right].
				$$
				Again, as the additional term is independent of $\mu$, we have $\partial_\mu G^M(\mu)(x)=\partial_\mu G(\mu;\mu_{\rm tar})(x)$.
				Consequently, we can also obtain that $\omega_M(R)=0$, and $\varepsilon_M(R)\le \frac{C}{M}$.

				\ms
				(3) Finally, we consider the so-called ``Regularized Wasserstein-type'' penalty functions. To be more precise, let us denote,  for $R>0$, $ \sB_R:=\big\{\mu\in\sP_2(\hR^d):  \int |x|^2\,\mu(dx)\le R\big\}$, and consider those penalty function $G$ satisfying Assumption \ref{assump1}-(ii) and that, for each $R>0$, there exists $C_R>0$, such that
				$$
				\left\{\ba{lll}
				|G(\mu;\nu)-G(\mu;\nu')| \le  C_R W_2(\nu,\nu'),   \qquad   \mu\in\sB_R, \ms\\
				\hE \Big[\Big| \partial_\mu G(\hP_X;\nu)(X) - \partial_\mu G(\hP_X;\nu')(X)\Big|^2\Big]\le  C_R W_2^2(\nu,\nu'), \q \mathbb E|X|^2\le R.
				\ea\right.
				$$
				Then, by  definition of $G^M$, and applying Jensen's inequality, one shows that
				$$
				\varepsilon_M(R)  \le  C_R\,\mathbb E W_2(\widehat\mu^M,\mu_{\rm tar})
				\q\mbox{and}\q
				\omega_M(R)  \le  C_R\,\mathbb E W_2^2(\widehat\mu^M,\mu_{\rm tar}).
				$$
				Furthermore, if $\mu_{\rm tar}$ has a finite $q$-th moment for some $q>4$, then the standard
				empirical Wasserstein estimates yield (cf. e.g., \cite{carmona2018probabilistic})
				$$
				\mathbb E W_2(\widehat\mu^M,\mu_{\rm tar})
				\le
				C
				\begin{cases}
					M^{-1/4}+M^{-(q-2)/(2q)}, & d<4,\\
					M^{-1/4}\sqrt{\log(1+M)}+M^{-(q-2)/(2q)}, & d=4,\\
					M^{-1/d}+M^{-(q-2)/(2q)}, & d>4.
				\end{cases}
				$$
				
				Hence $\varepsilon_M(R)$ and $\omega_M(R)$ satisfy Assumption \ref{assum:data_sample}-(iii) with the explicit bounds given as above. 
				\qed
			\end{Example}

			We are now ready to prove the convergence of the optimal control in Problem \ref{prob3} for each fixed $k$, as $M\rightarrow\infty$.
			To begin with, let us denote, for given $k\in\hN$, the optimal control for Problem \ref{softSBP} by $\widehat{\a}^k$, and that of Problem \ref{prob3} by $\widehat{\a}^{k,M}$, $M>0$. Also, in light of our previous discussion for Problem \ref{softSBP}, the optimal control for Problem \ref{prob3} should take the form 
			$\widehat{\alpha}^{k,M}_t = -\sigma(t)Y^{k,M}_t$, where $Y^{k,M}$ is the solution to the FBSDE \reff{eq:mkvfbsde_relaxed_km}. We have the following result. 
			\begin{Proposition}
				\label{prop:sample_approx} 
				Assume that Assumptions \ref{assump1} and  \ref{assum:data_sample} are in force. Then, for each $k\in\hN$, there exist constants $R_k, C_k>0$, independent of $M$, such that
				\bea
				\label{akMconv}
				\mathbb E\Big[\sup_{0\le t\le T}\Big|\widehat\alpha_t^{k,M}-\widehat\alpha_t^k\Big|^2\Big]+\mathbb E\int_0^T
				\Big|\widehat\alpha_t^{k,M}-\widehat\alpha_t^k\Big|^2 \d t \le C_k k^2\omega_M(R_k)\to 0, \q\mbox{as} ~M\to \infty.
				\eea
				
				Furthermore, it holds that $|V^{k,M}-V^k|\le k\varepsilon_M(R_k) \to 0$ as $
				M\to\infty$.
			\end{Proposition}
			The proof relies on Theorem~\ref{thm:global-fbsde-wellposedness}, and can be found in Appendix~\ref{sec:appendix}.

			\subsection{Finite particle system}
			\label{sec:particle_approximation}

			Fix the empirical measure $\widehat{\mu}^M$. Since the law $\hP_{X_T^\alpha}$ is not directly accessible, we approximate it by the empirical distribution of the terminal states of $N$ agents. All agents use the same feedback rule, evaluated at their respective states, and the objective is to minimize their average cost. This yields the following mean-field-type $N$-agent approximation of Problem \ref{prob3}, which is recovered as $N\to\infty$:

			\begin{Problem}
				\label{problem 4}
				Given a smooth penalty function $G$ and a constant $k>0$, minimize the cost functional:
				\begin{equation}
					\label{eq:relaxed_SB-objective1}
					J^{k,M,N}(\a)= \frac{1}{N}\sum_{j=1}^N \hE \left[ \frac{1}{2} \int_0^T |\a(s,X^{\a,j}_s)|^2 \d s  + k G^M(\mu_{T}^{N,\a}))\right]
				\end{equation}
				over all the closed-loop  strategy $\a= \a(t,\cdot)$ under the dynamics:
				\begin{equation}\label{eq:relaxed_SB-dynamics_finite_particle}
					\d X^{\a,j}_t = \big(b(t,X^{\a,j}_t) + \sigma(t) \a(t,X^{\a,j}_t)\big)\d t + \sigma(t)\d W^j_t, \quad  j=1,2,\cdots,N,
				\end{equation}
				where $\{W^j\}_{j=1,2,\cdots,N}$ are independent $d$-dimensional Brownian motions, $\{X^j_0\}_{j=1,2,\cdots,N}$ are i.i.d. with common distribution $\mu_{\rm ini}$ and are independent of the Brownian motions $\{W^j\}$, and $\mu_{T}^{N,\a}:= \frac{1}{N}\sum_{j=1}^N \delta_{X^{\a,j}_T}$ is the empirical distribution of the states of the $N$ agents at time $T$.
				\qed
			\end{Problem}
			
			We have the following result, which shows that $J^{k, M}$ can be approximated by $J^{k,M,N}$'s, and the error can be  explicitly expressed in terms of the number of particles $N$, the finite $q$-th moments of controlled dynamics $X_T^\alpha$, as well as the dimension $d$, thanks to the empirical Wasserstein estimates (cf. e.g., \cite{carmona2018probabilistic}) as we saw in 
			Example \ref{rem:GM_rates}-(3).
			
			\begin{Proposition}
				\label{prop:particle_approximation}
				Assume that Assumptions \ref{assump1},  and  \ref{assum:data_sample} are in force. Assume further that the     admissible control $\a$ is such that the corresponding controlled dynamics $X^\a$ satisfies
				\begin{equation}
					\label{eq:q_moment_assumption} 
					\E\big[|X_T^{\a}|^q\big]\;<\;\infty, \q \mbox{for some $q>2$.}
				\end{equation}
				Then there exists a constant $C=C(d,q)>0$ such that, for all $N\ge 2$,
				\begin{equation}\label{eq:quant_rate_cost_difference}
					\big|J_{k,M,N}(\alpha)-J_{k,M}(\alpha)\big|
					\;\le\;
					k\,L_0\, C \,\big(\E|X_T^\alpha|^q\big)^{1/q}\;\varepsilon_{d,q}(N),
				\end{equation}
				where $L_0$  is the $W_2$--Lipschitz constant of $G^M$ in Assumption~\ref{assum:data_sample}-(ii), and
				\begin{equation}\label{eq:eps_dq}
					\varepsilon_{d,q}(N):=
					\begin{cases}
						N^{-1/4}+N^{-(q-2)/(2q)}, & d<4,\\[3pt]
						N^{-1/4}\sqrt{\log N}+N^{-(q-2)/(2q)}, & d=4,\\[3pt]
						N^{-1/d}+N^{-(q-2)/(2q)}, & d>4.
					\end{cases}
				\end{equation}
				In particular, denoting  $\widehat\a^{M}$ to be the optimal strategy for Problem~\ref{prob3} and 
				$\widehat\a^{M,N}$ to be that of Problem~\ref{problem 4} we have, for all $N\ge 2$,
				\begin{equation}
					\label{eq:quant_rate_opt_gap}
					J_{k,M,N}(\widehat\a^{N,M})\;\le\;J_{k,M,N}(\widehat\a^{M})
					\;\le\;\inf_{\alpha\in \mathcal A}J_{k,M}(\alpha)
					\;+\;k\,L_0\, C \,\big(\E|X_T^{\h\a^{M}}|^q\big)^{1/q}\;\varepsilon_{d,q}(N),
				\end{equation}
				provided \eqref{eq:q_moment_assumption} holds for $\alpha=\widehat\a^{M}$.
				
			\end{Proposition}
			The proof uses the convergence rate of empirical measure in Wasserstein distance by~\cite{fournier2015rate}, and can be found in Appendix~\ref{sec:appendix}.

			
			

			\section{Numerical Algorithms}
			\label{sec:num_alg}
			In this section, we propose two algorithms for the numerical solution of Problem \ref{softSBP} in large-scale and high dimensions. For specific task such as image generation, which is of high-dimensionality, dimension reduction is usually applied first. Hence the Schr{\"o}dinger bridge is computed on the latent space, which is of smaller dimension.

			\subsection{MKV FBSDE Algorithm}
			\label{alg:FBSDE}

			\paragraph{New optimization problem for implementation.} The trick of the numerical approach to solve the FBSDE is to transform the forward-backward system into a forward system and an optimization problem aiming to satisfy the terminal condition required by the backward equation through the loss function  \cite{carmona2022convergence}
			$\hE[ |Y_T - k \pa_\mu G(\mathbb{P}_{X_T})(X_T)|^2]$,
			and the optimization problem is to search for the pair of functions $(y_0, z)$, $y_0 :\hR^d\to \hR^d$, $z:[0,T]\times \hR^d \to \hR^{d\times d}$,  which minimizes the cost functional
			\begin{equation}
				\label{eq:numerical_objective}
				J(y_0,z):=\hE\left[ \left|Y^{(y_0,z)}_T - k \pa_\mu G(\mathbb{P}_{X^{(y_0,z)}_T})(X^{(y_0,z)}_T)\right|^2\right],
			\end{equation}
			where the $(X,Y)$ solves
			\begin{align*}
				\label{eq:mkvfbsde_relaxed2}
				\begin{cases}
					\dd X^{(y_0,z)}_t = \left((b(t,X^{(y_0,z)}_t)- \sigma^2(t)Y^{(y_0,z)}_t\right)\dd t + \sigma(t) \dd W_t, \medskip\\
					\dd Y^{(y_0,z)}_t = - \nabla_x b(t,X^{(y_0,z)}_t)\cdot Y^{y_0,z}_t \dd t+ z(t,X^{(y_0,z)}_t) \dd W_t,\medskip\\
					X^{(y_0,z)}_0 \sim \mu_{\rm ini}, \quad Y^{(y_0,z)}_0 = y_0(X^{(y_0,z)}_0).
				\end{cases}
				\quad t\in[0,T],
			\end{align*}
			As implied by the above discussion, we will consider the uniform time grid and we will  approximate the functions $(y_0,z)$ by neural networks     $\mathcal{Y}_{\theta^y}:\mathbb{R}^d\rightarrow \mathbb{R}^d$ and $\mathcal{Z}_{\theta^z}:\mathbb{R}^{d+1}\rightarrow \mathbb{R}^d$ parameterized by $\theta^y$ and $\theta^z$, respectively. We assume the architecture for both networks are fixed, then we are aiming to minimize the cost functional over all the parameters $\theta:=(\theta^y,\theta^z)\in\Theta$.

			\paragraph{Time discretization.} We will consider a uniform grid $t_0=0<t_1<\cdots< t_N = T$ where $N$ denotes  the total number of partitions of the time interval. Denote each time increment $\Delta t= T/N$, $t_i = i \Delta t$, for $i=0,1,\dots, N$.

			Moreover, we shall approximate the distribution $\mathbb{P}_{X_T}$ in the cost functional by the empirical distribution of $B$ agents, where $B$ is often called the batch size. The empirical loss function over a batch is defined as
			\begin{equation}
				\label{eq:numerical_objective1}
				\widehat{\mathcal J}_B(\theta)= \left[
				\frac{1}{B} \sum_{j=1}^B | Y^j_{t_N}(\theta) - k\pa_\mu G(\bar\mu_{t_N}))(X^j_{t_N}(\theta))|^2
				\right],
			\end{equation}
			where the dynamic constraint for the $j$-th agent is given by 
			\eqref{eq:numerical_objective1}
			\begin{eqnarray}
				\label{eq:mkvfbsde_relaxed2}
				\begin{cases}
					X^j_{t_{i+1}} =  X^j_{t_i} + \left((b(t_i,X^j_{t_i})- \sigma(t_i)\sigma(t_i)^TY^j_{t_i}\right)\Delta t + \sigma(t_i)\Delta W^j_i, \medskip\\
					Y^j_{t_{i+1}} = Y^j_{t_i}  - \nabla_x b(t,X^j_{t_i})\cdot Y^j_{t_i} \Delta t + \mathcal{Z}_{\theta^z}(t,X^j_{t_i})\Delta W^j_i,
				\end{cases}
			\end{eqnarray}
			for all $i \in\{0,1,\dots,N\}$, and $\bar\mu_{t_N} = \frac{1}{B} \sum_{j=1}^B \delta_{X^j_{t_N}} $. 
			
			\paragraph{Neural network training.} 
			Note that we have two neural networks, one to parameterize $Y_0$ and the other to parameterize $Z$, with parameters summarized as $\theta:=(\theta^y,\theta^z)\in\Theta$.
			We then perform gradient descent update using the empirical loss function $\theta \leftarrow \theta - \eta_m \nabla_\theta \widehat{\mathcal J}_B(\theta)$
			with learning rate $\eta_m$. The full algorithm is summarized in Algorithm \ref{alg:fbsde_nn}.

			\subsection{MKV Control Algorithm}
			\label{alg:mkvc_nn}
			In addition to solving the McKean--Vlasov (MKV) control problem through the associated MKV FBSDE
			\eqref{eq:mkvfbsde_relaxed2}, we also implement a \emph{direct} numerical algorithm at the level of the
			MKV control problem, proposed in~\cite{soner2025learning}. Throughout, we work with the MKV control formulation 
			and construct a closed-loop feedback control using neural networks. The key idea is to approximate the
			flow of marginal laws by an interacting particle system and optimize the control parameters by minimizing
			an empirical approximation of the objective functional.
			\paragraph{MKV control problem and feedback form.}
			Recall the controlled diffusion dynamic~\eqref{eq:controlled-dynamics-follmer}.
			The relaxed Schr\"odinger bridge objective \eqref{eq:relaxed_SB-objective} can be seen as a special case
			of a terminal-cost MKV control problem. In the present section we describe the numerical method for the
			general form
			\begin{equation}\label{eq:mkvc_objective_general}
				\inf_{\alpha} \;
				\cJ(\alpha)
				:=
				\hE\Bigg[
				\int_0^T \Big(\tfrac12 |\alpha(t,X_t,\mu_t)|^2 + \ell(t,X_t,\mu_t)\Big)\,\dd t
				\;+\; k\,G(\mu_T)
				\Bigg],
			\end{equation}
			where $\ell$ is an optional running cost and $G$ is the terminal functional
			already appearing in \eqref{RSB_objective}. We consider feedback controls
			$\alpha:[0,T]\times \hR^d\times \sP_2(\hR^d)\to \hR^d$ and approximate $\mu_t$ by an empirical measure.
			
			\vspace{-5pt}
			\paragraph{Particle approximation and time discretization of the MKV dynamics.}
			Recall from Section~\ref{alg:FBSDE} that $B$ denotes the number of agents (i.e. batch size). Given a feedback control $\alpha$, the law $\mu_t$
			is approximated by the empirical measure
			$
			\bar\mu^B_t := \frac{1}{B}\sum_{j=1}^B \delta_{X^j_t}$,
			where $\{X^j_t\}_{j=1}^B$ solve the interacting system
			\begin{equation}\label{eq:particle_system}
				\dd X^j_t
				=
				b(t,X^j_t)\,\dd t
				+\sigma(t)\,\alpha(t,X^j_t,\bar\mu^B_t)\,\dd t
				+\sigma(t)\,\dd W^j_t,
				\qquad X^j_0 \stackrel{i.i.d.}{\sim}\mu_{\rm ini},
			\end{equation}
			with independent Brownian motions $\{W^j\}_{j=1}^B$. For numerical computation, we further use the uniform grid $t_0=0<t_1<\cdots<t_N=T$ with $\Delta t = T/N$ and apply the Euler-Maruyama scheme. Denote $\Delta W^j_i := W^j_{t_{i+1}}-W^j_{t_i} \sim \cN(0,\Delta t)$, independent over $i$ and $j$.
			Then \eqref{eq:particle_system} becomes
			\begin{equation}\label{eq:euler_particles}
				X^j_{t_{i+1}}
				=
				X^j_{t_i}
				+
				b(t_i,X^j_{t_i})\,\Delta t
				+
				\sigma(t_i)\,\alpha(t_i,X^j_{t_i},\bar\mu^B_{t_i})\,\Delta t
				+
				\sigma(t_i)\,\Delta W^j_i,
				\qquad i=0,\dots,N-1.
			\end{equation}
			
			\paragraph{Neural-network parametrization of the feedback control.}
			We approximate the feedback map $\alpha$ by a neural network with parameters $\theta$.
			Since $\mu$ is a probability measure, we do not feed $\mu$ directly into the network; instead we use a
			finite-dimensional \emph{feature map} of the empirical measure. 
			
			\smallskip
			Let
			$\Phi(\bar\mu^B_{t_i}) \in \hR^m$
			be a $m$-dimensional vector of statistics(e.g. moments, Fourier sequences, or a learned embedding). We then define the parameterized feedback control
			\begin{equation}\label{eq:alpha_nn}
				\alpha_\theta(t_i, X^j_{t_i}, \bar\mu^B_{t_i})
				:=
				\mathcal{A}_\theta\!\big(t_i,\,X^j_{t_i},\,\Phi(\bar\mu^B_{t_i})\big),
			\end{equation}
			where $\mathcal{A}_\theta $ is a standard feedforward network  with trainable parameters $\theta\in\Theta$.
			\paragraph{Empirical objective and optimization problem.}
			Given $\theta$, we simulate the particle system \eqref{eq:euler_particles} with feedback
			\eqref{eq:alpha_nn} and compute an empirical approximation of \eqref{eq:mkvc_objective_general}:
			\begin{equation}\label{eq:empirical_loss}
				\widehat{\cJ}_B(\theta)
				:=
				\frac1B\sum_{j=1}^B
				\Bigg[
				\sum_{i=0}^{N-1}
				\Big(\tfrac12\big|\alpha_\theta(t_i,X^j_{t_i},\bar\mu^B_{t_i})\big|^2
				+ f(t_i,X^j_{t_i},\bar\mu^B_{t_i})\Big)\Delta t
				\Bigg]
				\;+\;
				k\,G(\bar\mu^B_{t_N}),
			\end{equation}
			where $\bar\mu^B_{t_N}=\frac1B\sum_{j=1}^B \delta_{X^j_{t_N}}$.
			The numerical MKV control algorithm consists in solving
			\begin{equation}\label{eq:mkvc_nn_optim}
				\inf_{\theta\in\Theta}\; \hE\big[\widehat{\cJ}_B(\theta)\big],
			\end{equation}
			where the expectation is with respect to the simulated Brownian increments and the random initialization
			of particles.

			\begin{algorithm}[H]
				{\footnotesize
					\caption{Neural FBSDE for the Soft-Constrained Schr\"odinger Problem}
					\label{alg:fbsde_nn}
					\KwIn{Penalty weight $k>0$; time horizon $T$; number of steps $N$ ($\Delta t=T/N$); batch size $B$; initial law $\mu_{\rm ini}$; target law $\mu_{\rm tar}$ (through $G$); drift $b$, diffusion $\sigma$;  neural networks     $\mathcal{Y}_{\theta^y}$ and $\mathcal{Z}_{\theta^z}$ parameterized by $\theta^y$ and $\theta^z$ respectively; optimizer (Adam); number of iterations $K$.}
					\KwOut{Trained parameters $\h\theta:=(\h\theta^y,\h\theta^z)$ and networks $\mathcal{Y}_{\h\theta^y}$ and $\mathcal{Z}_{\h\theta^z}$.}
					
					\BlankLine
					Initialize $\theta:=(\theta^y,\theta^z) \leftarrow \theta^{(0)}$\;
					\For{$m=1,2,\dots,K$}{
						\tcp{(1) Sample initial particles and Brownian increments}
						Sample i.i.d.\ $X^{1}_{t_0},\dots,X^{B}_{t_0}\sim \mu_{\rm ini}$ and i.i.d. $\Delta W^{j}_i \sim \mathcal N(0,\Delta t)$ for $j=1,\dots,B$, $i=0,\dots,N-1$\;
						
						\tcp{(2) Forward simulate the interacting particle system}
						\For{$i=0,1,\dots,N-1$}{
							Set empirical measure $\bar\mu_{t_i} \leftarrow \frac1B\sum_{j=1}^B \delta_{X^j_{t_i}}$. Compute summary features $s_i \leftarrow \Phi(\bar\mu_{t_i})$\;
							\For{$j=1,\dots,B$}{
								Set $Y_0^j=\mathcal{Y}_{\theta^y}(X^j_0))$. Update $X^j_{t_{i+1}}$ and  $Y^j_{t_{i+1}}$  by~\eqref{eq:mkvfbsde_relaxed2};
							}
						}
						Set $\bar\mu_T \leftarrow \frac1B\sum_{j=1}^B \delta_{X^j_{t_N}}$\;
						
						\tcp{(3) Evaluate the empirical loss}
						Compute $\widehat{\mathcal J}_B(\theta)$ by~\eqref{eq:numerical_objective1}\;
						\tcp{(4) Gradient step}
						Update parameters $\theta \leftarrow \theta - \eta_m \nabla_\theta \widehat{\mathcal J}_B(\theta)$ using automatic differentiation and Adam\;
					}
					\Return{$\h\theta \leftarrow \theta$ }\;
				}
			\end{algorithm}
			
			\begin{Remark}[Choice of measure features]
				The feature map $\Phi(\bar\mu^B_{t})$ can be chosen in several ways. Common options include: (i) finitely
				many moments $(\int x^p\,d\bar\mu^B_t)_{p\le P}$;
				(ii) evaluations against a separating class of test functions $(\int \varphi_\ell\,d\bar\mu^B_t)_{\ell=1}^m$;
				or (iii) a learned embedding (e.g.\ a small ``statistics network'').
				The choice of $\Phi$ trades off expressiveness and stability.
			\end{Remark}

			\begin{algorithm}[H]
				{\footnotesize \caption{Neural MKV Control for the Soft-Constrained Schr\"odinger Problem}
					\KwIn{Penalty weight $k>0$; time horizon $T$; number of steps $N$ ($\Delta t=T/N$); batch size $B$; initial law $\mu_{\rm ini}$; target law $\mu_{\rm tar}$ (through $G$); drift $b$, diffusion $\sigma$; feature map $\Phi$; NN policy class $\alpha_\theta(t,x,\Phi(\bar\mu))$; optimizer (Adam); number of iterations $K$.}
					\KwOut{Trained parameters $\h\theta$ and feedback control $\alpha_{\h\theta}$.}
					
					\BlankLine
					Initialize $\theta \leftarrow \theta^{(0)}$\;
					\For{$m=1,2,\dots,K$}{
						\tcp{(1) Sample initial particles and Brownian increments}
						Sample i.i.d.\ $X^{1}_{t_0},\dots,X^{B}_{t_0}\sim \mu_{\rm ini}$\;
						Sample i.i.d.\ Brownian increments $\Delta W^{j}_i \sim \mathcal N(0,\Delta t)$ for $j=1,\dots,B$, $i=0,\dots,N-1$\;
						
						\tcp{(2) Forward simulate the interacting particle system}
						\For{$i=0,1,\dots,N-1$}{
							Set empirical measure $\bar\mu_{t_i} \leftarrow \frac1B\sum_{j=1}^B \delta_{X^j_{t_i}}$. Compute summary features $s_i \leftarrow \Phi(\bar\mu_{t_i})$\;
							\For{$j=1,\dots,B$}{
								Set $\alpha^j_i \leftarrow \alpha_\theta(t_i, X^j_{t_i}, s_i)$. Update $X^j_{t_{i+1}}$ by~\eqref{eq:euler_particles};
							}
						}
						Set $\bar\mu_T \leftarrow \frac1B\sum_{j=1}^B \delta_{X^j_{t_N}}$\;
						
						\tcp{(3) Evaluate the empirical loss}
						Compute $ \widehat{\mathcal J}_B(\theta)$ by~\eqref{eq:empirical_loss}\;
						\tcp{(4) Gradient step}
						Update parameters $\theta \leftarrow \theta - \eta_m \nabla_\theta \widehat{\mathcal J}_B(\theta)$ using automatic differentiation and Adam\;
					}
					\Return{$\h\theta \leftarrow \theta$ and control $\alpha_{\h\theta}$}\;
				}
			\end{algorithm}

			\section{Experiments}
			\subsection{Unconditional generation on MNIST}
			We demonstrate the soft-constrained Schr\"odinger formulation on the MNIST dataset\cite{lecun2010mnist} by learning a controlled diffusion that transports an easy-to-sample prior (standard Gaussian) to the MNIST data distribution. To keep the state dimension manageable, we work in a learned latent space: a Variational Auto-Encoder(VAE) first compresses images into a $16$-dimensional representation, and then we solve the mean-field control (MFC) problem in latent space, with the two numerical algorithms in Section~\ref{sec:num_alg}.
			
			We take the initial distribution $\mu_{\rm ini}=\cN(0,I_{16})$, $\sigma=1$, $T=1$, $\Delta t=1/40$, $G(\mu):=\rho_*^2(\mu,\h{\mu}_{\rm tar})$ as in Example~\ref{ex: fourier-wasserstein} with terminal penalty parameter $k=20$. For the neural networks, we use two-hidden-layer Multi-Layer Perceptrons(MLP), where each hidden layer contains 256 neurons. For the training process, we apply Adam algorithm with $20$ epochs and learning rate $10^{-3}$. For both algorithms, we simulate $B=128$ particles for dynamics approximation.
			
			After training, the new generated results are shown in Figure~\ref{fig:MNIST}.
			\begin{figure}
				\centering
				\includegraphics[width=0.4\linewidth]{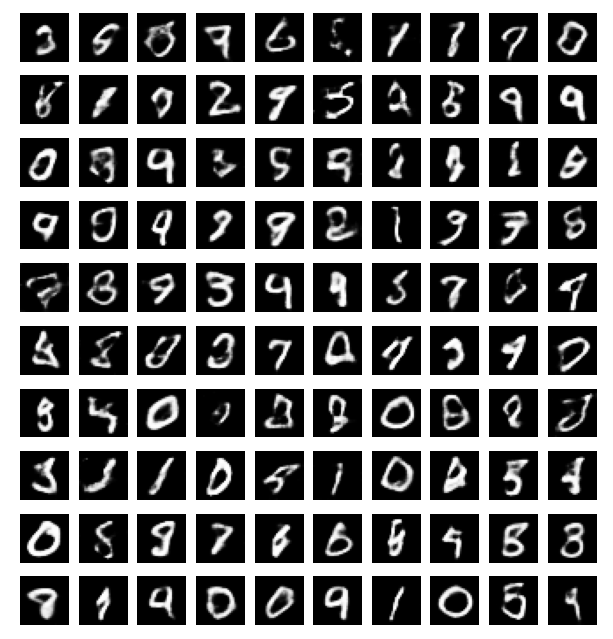}
				\includegraphics[width=0.32\linewidth]{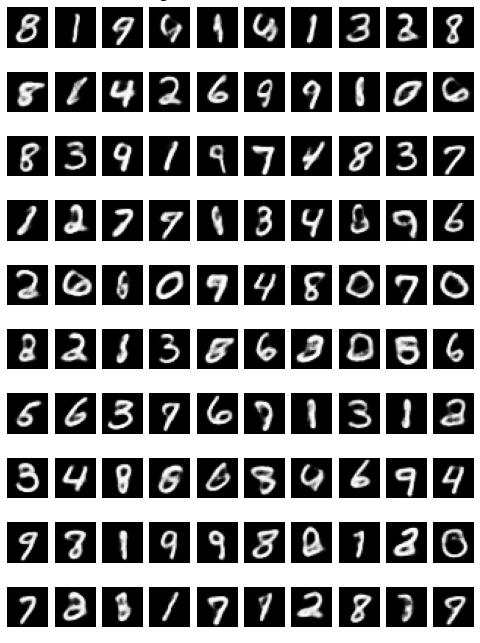}
				\caption{ Left: FBSDE.  Right: Direct MFC.}
				\label{fig:MNIST}
			\end{figure}

			\subsection{Adult-face to children-face generation}
			Following \cite{alouadi2026lightsbbmbridgingschrodingerbass}, we formulate the adult-face to children-face \emph{distribution-to-distribution} image translation task
			through our \emph{soft-constrained Schr\"odinger problem}. In this formulation, there is no identity-preservation label in the dataset, and the goal is not paired
			translation, but rather to learn a population-level control
			that steers the law of a stochastic process from the adult-face distribution to the children-face
			distribution. After training, we apply the learned optimal control to \emph{new} testing adult-face inputs and
			simulate the controlled dynamics via Euler--Maruyama to generate visually plausible children faces.
			
			\paragraph{Configurations}
			The training dataset consists of an \emph{adult-face} dataset (age $\ge\, 18$) $\cD_{\rm ini}$ and a \emph{children-face} dataset (age $<\, 18$) $\cD_{\rm tar}$, each contains $n_{\rm ini}=n_{\rm tar}=60000$ images. All images are aligned/cropped and resized to a common resolution of $d_{\rm real}=1024\times1024 $. We first apply a fixed dimension-reduction map $\Psi$ to embed faces into a $d_{\rm lat}=512$ dimensional latent space using Variational Autoencoder~\cite{VAE}. Hence, the dimension-reduction map $\Psi: \hR^{d_{\rm real}}\to \hR^{d_{\rm lat}}$. We denote the induced empirical distributions in $\R^{d_{\rm lat}}$ 
			$ 
			\mu_{\rm ini}\approx \frac1{n_{\rm ini}}\sum_{i=1}^{n_{\rm ini}}\delta_{\Psi(I^{\rm adult}_i)}
			$ and $
			\mu_{\rm tar}\approx \frac1{n_{\rm tar}}\sum_{i=1}^{n_{\rm tar}}\delta_{\Psi(I^{\rm child}_i)},
			$
			where for every $i\in \{1,\ldots,n_{\rm ini}\}$, $I_i^{\rm adult}\in \cD_{\rm ini},\,I_i^{\rm child}\in\cD_{\rm tar}$.  All Schr\"odinger/MKV-control computations are performed in this reduced space. We define the penalty function 
			$$
			G(\mu)=SW_2(\mu,\mu_{\rm tar})^2\quad \text{or} \quad G(\mu)= \rho_{\lambda}(\mu,\mu_{\rm tar})^2,
			$$
			where $SW_2$ denotes the sliced Wasserstein-two distance, and $\rho_{\lambda}$ is defined in Example~\ref{ex: fourier-wasserstein}. For the soft-constrained Schr\"odinger problem, we set $T=1$ with $50$-step time grid and set $\sigma=0.05$, $k=50$.
			
			\smallskip
			To parametrize the feedback control, we use feedforward neural networks with three hidden layers with SiLU activation function, each with $1024$ hidden neurons. The training uses Adam gradient descent algorithm, takes $500\sim 1000$ epochs with learning rate $2e-4$. 
			
			\paragraph{Numerical results}
			After training, we take $N_{\rm test}$ initial configurations from the testing dataset(adult faces), apply Euler-Maruyama scheme to simulate the MKV-dynamics~\eqref{eq:particle_system}. We illustrate the results on $N_{\rm test}=36$ testing images. 
			For FBSDE algorithm~\ref{alg:FBSDE}, we illustrate the results on $N_{\rm test}=36$ images in Figure~\ref{fig:oldandyound_FBSDE}. For MKV algorithm~\ref{alg:mkvc_nn}, the results are illustrated in Figure~\ref{fig:oldandyoung}. 
			\begin{figure}[H]
				\centering
				\includegraphics[width=0.35\linewidth]{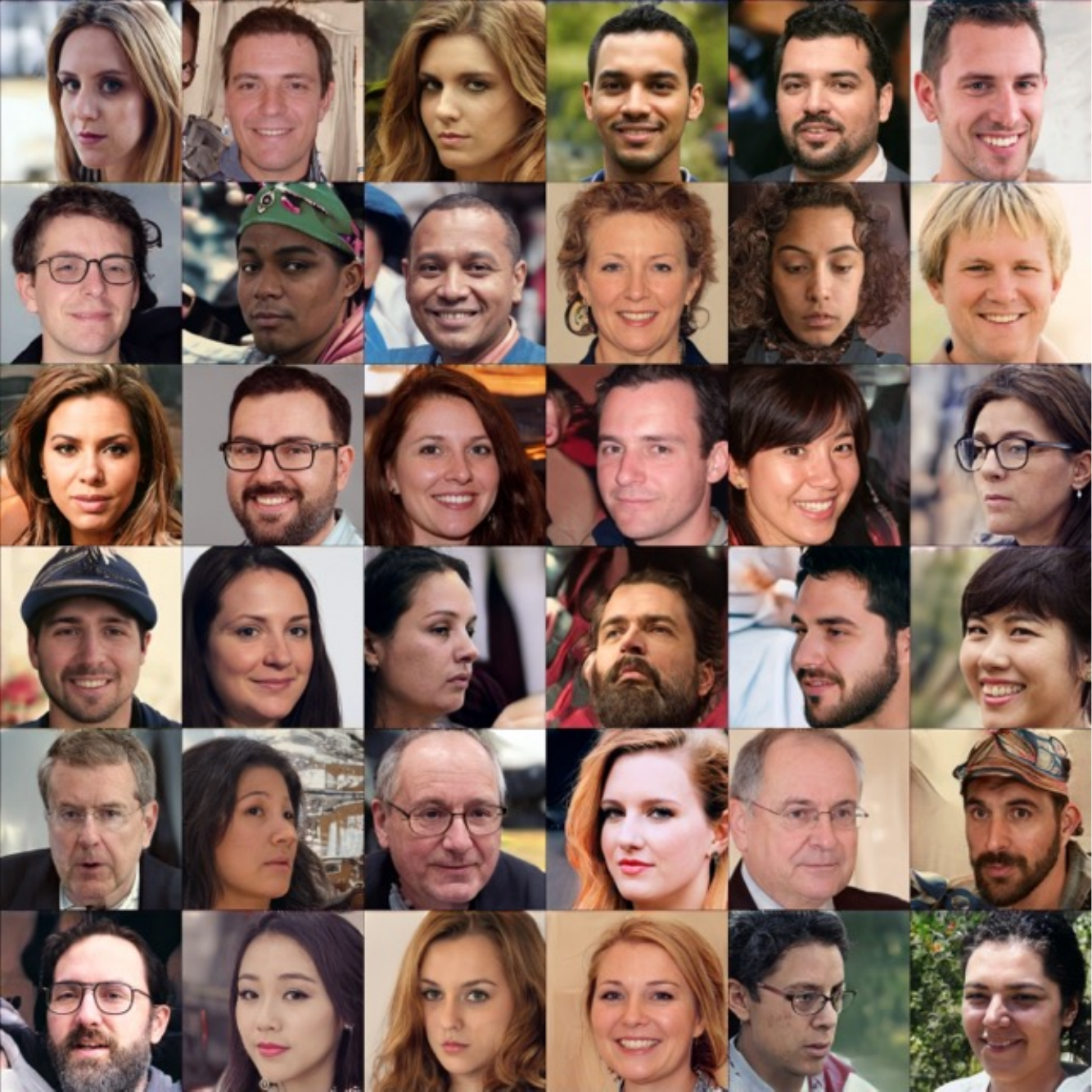}
				\includegraphics[width=0.35\linewidth]{Figures/FBSDE_Child.pdf}
				\caption{FBSDE. Left: Initial images(Adult faces).   Right: Terminal images(Child faces).
				}
				\label{fig:oldandyound_FBSDE}
			\end{figure}
			\begin{figure}[H]
				\centering
				\includegraphics[width=0.35\linewidth]{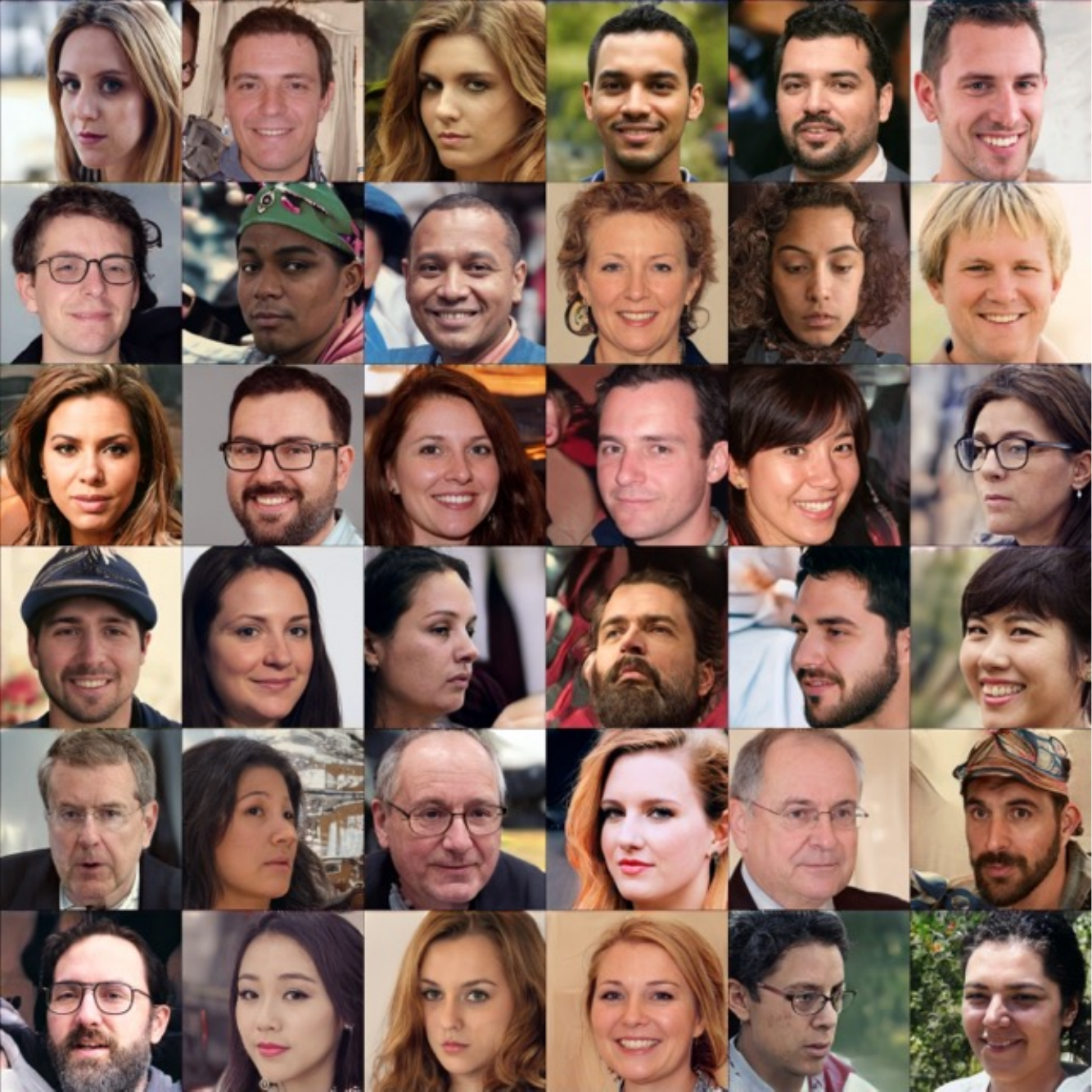}
				\includegraphics[width=0.35\linewidth]{Figures/MFC_Child.pdf}
				\caption{MKV. Left: Initial images(Adult faces).   Right: Terminal images(Child faces).
				}
				\label{fig:oldandyoung}
			\end{figure}

			{\footnotesize
				\bibliographystyle{abbrv}
				\bibliography{references}
			}

			\newpage
			
			\appendix
			\section{Additional Proofs}
			\label{sec:appendix}
			\begin{proof}[Proof of~\ref{lem:path-law-existence-uniqueness}]
				Let $\{\hQ_n\}_{n\geq1}\subset\sQ_{\mu_{\mathrm{ini}}}$ be a minimizing sequence of \eqref{eq:path-law-objective}. Since $G\geq0$, we can assume that $D_{\mathrm{KL}}(\hQ_n\Vert \hP)\le V^k+1$, $n\ge 1$. Since $\Omega$ is Polish, the level set of relative-entropy $D(\cd\|\hP)$ is weakly compact (see, e.g., \cite{DupuisEllis1997}). Thus, possibly along a subsequence and may assume itself, $\hQ_n$ converges weakly to some $\hQ\in\sP(\Omega)$.
				
				Next, note that the projection mapping $\omega\mapsto\omega(0)$ is continuous, and the mapping $\omega\mapsto\sup_{0\leq t\leq T}|\omega_t|^2$ is lower semicontinuous, the Portmanteau theorem implies that the limit $\hQ$ preserves the properties: 
				$ \hQ_{X_0}=\mu_{\mathrm{ini}}$ and $\mathbb E^{\hQ}\left[\|X\|^2_{\hC([0,T])}\right]<\infty$.
				That is, $\hQ\in\sQ_{\mu_{\mathrm{ini}}}$. 
				
				Finally, note that  the relative entropy  $D(\,\cd\,\|\hP)$
				is weakly lower semicontinuous and so is $G$, by Assumption~\ref{assump1}, and the projection mapping $\hQ\mapsto \hQ_{X_T}$ is weakly continuous,
				we thus have
				$$
				J^k(\hQ) \leq \liminf_{n\to\infty}J^k(\hQ_n)=V^k.
				$$
				In other words,  $\hQ$ is a minimizer of \eqref{eq:path-law-objective}.
				
				To prove uniqueness, we need only note that 
				the set $\mathcal Q_{\mu_{\mathrm{ini}}}$ is convex, and
				the map $\hQ\mapsto D_{\mathrm{KL}}(\hQ\Vert\hP)$ is strictly convex  on $\mathcal Q_{\mu_{\mathrm{ini}}}$, whereas $G$ is convex  
				(under linear interpolation of probability measures). Thus,  $J^k(\cd)$ is strictly convex on $\mathcal Q_{\mu_{\mathrm{ini}}}$, which implies the uniqueness of the minimizer. 
			\end{proof}

			\begin{proof}[Proof of Proposition \ref{prop:sample_approx}.]
				First note that under the assumptions, both FBSDE  \reff{eq:mkvfbsde_relaxed_km} is
				well-posed, thanks to Theorem~\ref{thm:global-fbsde-wellposedness}. Furthermore,  since $\partial_\mu G$ and $\partial_\mu G^M$ are uniformly bounded, uniformly in $M$, we have
				$$
				|Y_T^k|=k|\partial_\mu G(P_{X_T^k})(X_T^k)|\le kL_0, \qq
				|Y_T^{k,M}|=k|\partial_\mu G^M(P_{X_T^{k,M}})(X_T^{k,M})|\le kL_0, \q M>0.
				$$
				The backward equations are linear: $dY_t=-A_tY_t\,dt+Z_t\,dW_t,
				$ with $\|A_t\|\le L$. The conditional representation for this linear BSDE and Assumption \ref{assum:data_sample}-(ii) give
				\[
				|Y_t|\le e^{L(T-t)}\mathbb E[|Y_T|\mid\mathcal F_t]
				\le kL_0e^{LT}.
				\]The same argument applies to \(Y^M\). Thus, there exists a constant $C_k$ independent of $M$, such that
				\bea
				\label{YCk}
				\|Y^k\|_{L^\infty([0,T]\times\Omega)}+\sup_{M\ge1}\|Y^{k,M}\|_{L^\infty([0,T]\times\Omega)}\le C_k.
				\eea
				Since $|b(t,x)|\le \sup_{s\le T}|b(s,0)|+L|x|$ by Assumption \ref{assump1}(1)
				and \(Y,Y^M,\sigma\) are uniformly bounded, the forward equations, the martingale maximal inequality, and Gronwall’s inequality yield
				\[
				\mathbb E\sup_{t\le T}|X^k_t|^2
				+\sup_{M\ge1}\mathbb E\sup_{t\le T}|X_t^{k,M}|^2
				\le C\bigl(1+\mathbb E|\xi_0|^2+k^2L_0^2e^{2LT}\bigr)
				=:R_k.
				\]
				Hence it holds that
				\bea
				\label{XRk}
				\mathbb E|X_T^k|^2+\sup_{M\ge1}\mathbb E|X_T^{k,M}|^2 \le R_k.
				\eea
				
				To see the difference between the solutions to the two FBSDEs. Let us denote
				$$
				I_M:=k\Big[\partial_\mu G^M(P_{X_T^{k,M}})(X_T^{k,M})-\partial_\mu G(P_{X_T^{k,M}})(X_T^{k,M})\Big],\quad 
				Y_T^{k,M}=k\partial_\mu G(P_{X_T^{k,M}})(X_T^{k,M})+I_M.
				$$
				Then by the definition of $\omega_M(R)$ (see \reff{eoMR}) and estimate (\ref{XRk}), we easily deduce that
				$
				\mathbb E|I_M|^2\le k^2\omega_M(R_k)$.
				
				By Itô's formula and the equation for $v^\mu$ in Lemma~\ref{lem:frozen-log-transform-regularity},
				\begin{equation}\label{eq:frozen-bsde}
					\begin{aligned}
						&d\left[Y_t^{k,M}
						-\nabla_xv^{\mathbb P_{X_T^{k,M}}}(t,X_t^{k,M})\right]\\
						&=-\left[(D_xb(t,X_t^{k,M}))^\top
						-\nabla_x^2v^{\mathbb P_{X_T^{k,M}}}(t,X_t^{k,M})a(t)\right]\cdot\left[Y_t^{k,M}
						-\nabla_xv^{\mathbb P_{X_T^{k,M}}}(t,X_t^{k,M})\right]dt\\
						&\quad+\left[Z_t^{k,M}
						-\nabla_x^2v^{\mathbb P_{X_T^{k,M}}}(t,X_t^{k,M})\sigma(t)\right]dW_t,
					\end{aligned}
				\end{equation}
				with terminal value $Y_T^{k,M}
				-\nabla_xv^{\mathbb P_{X_T^{k,M}}}(T,X_T^{k,M})=I_M.$
				The coefficient in \eqref{eq:frozen-bsde} is bounded independently
				of $M$ by Assumption~\ref{assump1}(1)(i) and \eqref{eq:uniform-frozen-field-bound}. The linear BSDE estimate
				and Girsanov's formula therefore give
				\begin{equation}\label{eq:frozen-entropy}
					\begin{aligned}
						D_{\mathrm{KL}}\left(
						\mathbb P_{X^{k,M}}
						\,\middle\|\,
						\widehat{\mathbb Q}^{\mathbb P_{X_T^{k,M}}}
						\right)=\frac12\mathbb E\int_0^T
						\left|\sigma(t)^\top\left[
						Y_t^{k,M}
						-\nabla_xv^{\mathbb P_{X_T^{k,M}}}(t,X_t^{k,M})
						\right]\right|^2dt
						\le C_k\mathbb E|I_M|^2.
					\end{aligned}
				\end{equation}
				Here $\mathbb P_{X^{k,M}}$ and $\mathbb P_{X^k}$ denote
				the laws of the entire paths.
				
				By Proposition~\ref{prop:self-consistent-gibbs}, Lemma~\ref{lem:conditional-gibbs}, and the definition of $B_G$,
				\begin{equation}\label{eq:entropy-comparison}
					\begin{aligned}
						D_{\mathrm{KL}}\left(
						\mathbb P_{X^{k,M}}\middle\|\mathbb P_{X^k}\right)&\le J^k(\widehat\alpha^{k,M})-J^k(\widehat\alpha^k)\\
						&=D_{\mathrm{KL}}\left(
						\mathbb P_{X^{k,M}}
						\,\middle\|\,
						\widehat{\mathbb Q}^{\mathbb P_{X_T^{k,M}}}
						\right)
						-D_{\mathrm{KL}}\left(
						\mathbb P_{X^k}
						\,\middle\|\,
						\widehat{\mathbb Q}^{\mathbb P_{X_T^{k,M}}}
						\right)
						-kB_G\left(\mathbb P_{X_T^k}\mid\mathbb P_{X_T^{k,M}}\right)\\
						&\le D_{\mathrm{KL}}\left(
						\mathbb P_{X^{k,M}}
						\,\middle\|\,
						\widehat{\mathbb Q}^{\mathbb P_{X_T^{k,M}}}
						\right)
						\le C_k\mathbb E|I_M|^2,
					\end{aligned}
				\end{equation}
				where the last two inequalities use $B_G\ge0$ and
				\eqref{eq:frozen-entropy}.
				
				The map
				$x\mapsto b(t,x)-a(t)\nabla_xv^{\mathbb P_{X_T^k}}(t,x)$
				is globally Lipschitz by \eqref{eq:uniform-frozen-field-bound}. Subtracting the forward equations
				with the common initial variable and Brownian motion, Gronwall's
				inequality and Girsanov's formula yield
				\begin{equation}\label{eq:state-comparison}
					\begin{aligned}
						\mathbb E\sup_{t\le T}|X_t^{k,M}-X_t^k|^2
						&\le C_k\mathbb E\int_0^T
						\left|\sigma(t)^\top\left[
						Y_t^{k,M}
						-\nabla_xv^{\mathbb P_{X_T^k}}(t,X_t^{k,M})
						\right]\right|^2dt\\
						&=2C_kD_{\mathrm{KL}}\left(
						\mathbb P_{X^{k,M}}\middle\|\mathbb P_{X^k}\right)
						\le C_k\mathbb E|I_M|^2.
					\end{aligned}
				\end{equation}
				
				The difference of the backward equations is linear, with bounded
				coefficient $(D_xb(t,X_t^{k,M}))^\top$. Its forcing and terminal
				value satisfy, by Assumptions~\ref{assump1}(1)(i), \ref{assum:data_sample}(ii), and \eqref{YCk},
				\[
				\begin{aligned}
					\left|\left[(D_xb(t,X_t^{k,M}))^\top
					-(D_xb(t,X_t^k))^\top\right]Y_t^k\right|
					&\le C_k|X_t^{k,M}-X_t^k|,\\
					\mathbb E|Y_T^{k,M}-Y_T^k|^2
					&\le 2k^2L_0\mathbb E|X_T^{k,M}-X_T^k|^2
					+2\mathbb E|I_M|^2.
				\end{aligned}
				\]
				The linear BSDE estimate and \eqref{eq:state-comparison}
				therefore imply
				\begin{equation}\label{eq:fbsde-stability}
					\mathbb E\Big[
					\sup_{0\le t\le T}|X_t^{k,M}-X_t^k|^2
					+\sup_{0\le t\le T}|Y_t^{k,M}-Y_t^k|^2
					+\int_0^T|Z_t^{k,M}-Z_t^k|^2\,dt
					\Big]
					\le C_k\mathbb E|I_M|^2
					\le C_k k^2\omega_M(R_k).
				\end{equation}

				Moreover, since $\widehat\alpha_t^{k,M}-\widehat\alpha_t^k=-\sigma(t)^\top (Y_t^{k,M}-Y_t^k )$, we have
				$
				\left|\widehat\alpha_t^{k,M}-\widehat\alpha_t^k\right|^2
				\le
				\|\sigma\|_\infty^2
				|Y_t^{k,M}-Y_t^k|^2$.
				Therefore,
				$$
				\mathbb E\Big[\sup_{0\le t\le T}\Big|\widehat\alpha_t^{k,M}-\widehat\alpha_t^k\Big|^2\Big]+
				\mathbb E\int_0^T
				|\widehat\alpha_t^{k,M}-\widehat\alpha_t^k |^2 dt\le C_k k^2\omega_M(R_k).
				$$
				
				It remains to prove the convergence of the optimal values. To this end, we first note that by using  the optimality of $\widehat\alpha^{k,M}$ (resp. $\widehat\alpha^k$) for $J^{k,M}$ (resp. $J^k$) we have
				\bea
				\label{VkM-k1}
				V^{k,M}-V^k=J^{k,M}(\widehat\alpha^{k,M})-J^k(\widehat\alpha^k)
				\le J^{k,M}(\widehat\alpha^k)-J^k(\widehat\alpha^k).
				\eea
				Next, recalling the definition of $J^k(\cd)$ (\ref{RSB_objective}) and $J^{k, M}(\cd)$ (\ref{eq:relaxed_SB-objective}), we see that, under the same control $\widehat\alpha^k$, it holds that
				\bea
				\label{JkM-k}
				J^{k,M}(\widehat\alpha^k)-J^k(\widehat\alpha^k)=k\big[G^M(\hP_{X_T^k})-G(\hP_{X_T^k})\big].
				\eea
				Using the definition of $\varepsilon_M(R)$ and the moment bound, we derive from (\ref{VkM-k1}) and (\ref{JkM-k}) that
				\bea
				\label{VkM-k2}
				V^{k,M}-V^k\le k\varepsilon_M(R_k).
				\eea
				By switching the roles of $V^{k, M}$ and $V^k$ we can similarly  argue that
				\bea
				\label{Vk-kM}
				V^k-V^{k,M}
				\le J^k(\widehat\alpha^{k,M})-J^{k,M}(\widehat\alpha^{k,M})=k\big[G(P_{X_T^{k,M}})-G^M(P_{X_T^{k,M}})\big]\le
				k\varepsilon_M(R_k). 
				\eea
				Combining (\ref{VkM-k2}) and (\ref{Vk-kM}), we obtain
				$
				|V^{k,M}-V^k|
				\le
				k\varepsilon_M(R_k)$.
				Finally, by Assumption~\ref{assum:data_sample}, we conclude that
				$V^{k,M}\to V^k$, as $M\to\infty$, proving the proposition.
			\end{proof}
			
			\begin{proof}[Proof of Proposition~\ref{prop:particle_approximation}]
				Fix $\alpha\in \mathcal A$.
				Because the state dynamics \eqref{eq:relaxed_SB-dynamics_finite_particle} are decoupled across $j$ for a given feedback $\alpha$,
				the terminal states $\{X_T^{\alpha,j}\}_{j=1}^N$ are i.i.d.\ with common law $\mu_T^\alpha$ and
				$\mu_T^{N,\alpha}$ is their empirical measure.
				
				By exchangeability,
				\begin{eqnarray*}
					&&   J_{k,M,N}(\alpha)
					=\E\Big[\frac12\int_0^T|\alpha(t,X_t^\alpha)|^2\,dt\Big]
					+k\,\E\big[G_M(\mu_T^{N,\alpha})\big],\\
					&&J_{k,M}(\alpha)
					=\E\Big[\frac12\int_0^T|\alpha(t,X_t^\alpha)|^2\,dt\Big]
					+k\,G_M(\mu_T^\alpha).
				\end{eqnarray*}
				Hence
				\[
				\big|J_{k,M,N}(\alpha)-J_{k,M}(\alpha)\big|
				= k\Big|\E\big[G_M(\mu_T^{N,\alpha})\big]-G_M(\mu_T^\alpha)\Big|
				\le k\,\E\big|G_M(\mu_T^{N,\alpha})-G_M(\mu_T^\alpha)\big|.
				\]
				Using the $W_2$--Lipschitz property of $G_M$ in Assumption~\ref{assum:data_sample}(ii), $\big|G_M(\mu_T^{N,\alpha})-G_M(\mu_T^\alpha)\big|
				\le L_0\,W_2(\mu_T^{N,\alpha},\mu_T^\alpha),$ thus
				\begin{equation}\label{eq:reduce_to_W2}
					\big|J_{k,M,N}(\alpha)-J_{k,M}(\alpha)\big|
					\le k\,L_0\,\E\Big[W_2(\mu_T^{N,\alpha},\mu_T^\alpha)\Big].
				\end{equation}
				
				To bound the expectation of the empirical Wasserstein distance, we invoke a standard
				non-asymptotic estimate for i.i.d.\ empirical measures in $W_2$ under a finite $q$-moment:
				by \cite[Theorem~1]{fournier2015rate}, assumption \eqref{eq:q_moment_assumption}
				implies that there exists $C=C(d,q)>0$ such that
				\[
				\E\Big[W_2(\mu_T^{N,\alpha},\mu_T^\alpha)\Big]
				\le C\,\big(\E|X_T^\alpha|^q\big)^{1/q}\;\varepsilon_{d,q}(N),
				\]
				with $\varepsilon_{d,q}(N)$ given by \eqref{eq:eps_dq}.
				Plugging this bound into \eqref{eq:reduce_to_W2} yields \eqref{eq:quant_rate_cost_difference}.
				
				Finally, if $\h \a_{k,M}$ is optimal for Problem~\ref{prob3}, then $J_{k,M}(\h \a^{k,M})=\inf_{\alpha\in \mathcal A}J_{k,M}(\alpha).$
				Moreover, by optimality of $\h\alpha^{M,N}$ for Problem~\ref{problem 4},
				$J_{k,M,N}(\h\alpha^{M,N})\le J_{k,M,N}(\h\alpha^{M})$ for each $N$.
				Combining these inequalities with \eqref{eq:quant_rate_cost_difference} applied to $\alpha=\h\alpha^{M}$
				gives \eqref{eq:quant_rate_opt_gap}.
			\end{proof}

		\end{document}